%%%%%%%%%%%%%%%%%%%%%%%%%%%%%%% including smallfonts.tex %%%%%%%%%%%%%%%%%%%%%%%%%%%%%%%
%smallfonts.tex
%
\newskip\ttglue
\font\fiverm=cmr5
\font\fivei=cmmi5
\font\fivesy=cmsy5
\font\fivebf=cmbx5
\font\sixrm=cmr6
\font\sixi=cmmi6
\font\sixsy=cmsy6
\font\sixbf=cmbx6
\font\sevenrm=cmr7
\font\eightrm=cmr8
\font\eighti=cmmi8
\font\eightsy=cmsy8
\font\eightit=cmti8
\font\eightsl=cmsl8
\font\eighttt=cmtt8
\font\eightbf=cmbx8
\font\ninerm=cmr9
\font\ninei=cmmi9
\font\ninesy=cmsy9
\font\nineit=cmti9
\font\ninesl=cmsl9
\font\ninett=cmtt9
\font\ninebf=cmbx9
\font\twelverm=cmr12
\font\twelvei=cmmi12
\font\twelvesy=cmsy12
\font\twelveit=cmti12
\font\twelvesl=cmsl12
\font\twelvett=cmtt12
\font\twelvebf=cmbx12

%% EIGHT POINT FONT FAMILY

\def\eightpoint{\def\rm{\fam0\eightrm}
  \textfont0=\eightrm \scriptfont0=\sixrm \scriptscriptfont0=\fiverm
  \textfont1=\eighti  \scriptfont1=\sixi  \scriptscriptfont1=\fivei
  \textfont2=\eightsy  \scriptfont2=\sixsy  \scriptscriptfont2=\fivesy
  \textfont3=\tenex  \scriptfont3=\tenex  \scriptscriptfont3=\tenex
  \textfont\itfam=\eightit  \def\it{\fam\itfam\eightit}
  \textfont\slfam=\eightsl  \def\sl{\fam\slfam\eightsl}
  \textfont\ttfam=\eighttt  \def\tt{\fam\ttfam\eighttt}
  \textfont\bffam=\eightbf  \scriptfont\bffam=\sixbf
    \scriptscriptfont\bffam=\fivebf  \def\bf{\fam\bffam\eightbf}
  \tt  \ttglue=.5em plus.25em minus.15em
  \normalbaselineskip=9pt
  \setbox\strutbox=\hbox{\vrule height7pt depth2pt width0pt}
  \let\sc=\sixrm  \let\big=\eightbig \normalbaselines\rm}

\def\eightbig#1{{\hbox{$\textfont0=\ninerm\textfont2=\ninesy
        \left#1\vbox to6.5pt{}\right.$}}}

%% NINE POINT FONT FAMILY

\def\ninepoint{\def\rm{\fam0\ninerm}
  \textfont0=\ninerm \scriptfont0=\sixrm \scriptscriptfont0=\fiverm
  \textfont1=\ninei  \scriptfont1=\sixi  \scriptscriptfont1=\fivei
  \textfont2=\ninesy  \scriptfont2=\sixsy  \scriptscriptfont2=\fivesy
  \textfont3=\tenex  \scriptfont3=\tenex  \scriptscriptfont3=\tenex
  \textfont\itfam=\nineit  \def\it{\fam\itfam\nineit}
  \textfont\slfam=\ninesl  \def\sl{\fam\slfam\ninesl}
  \textfont\ttfam=\ninett  \def\tt{\fam\ttfam\ninett}
  \textfont\bffam=\ninebf  \scriptfont\bffam=\sixbf
    \scriptscriptfont\bffam=\fivebf  \def\bf{\fam\bffam\ninebf}
  \tt  \ttglue=.5em plus.25em minus.15em
  \normalbaselineskip=11pt
  \setbox\strutbox=\hbox{\vrule height8pt depth3pt width0pt}
  \let\sc=\sevenrm  \let\big=\ninebig \normalbaselines\rm}

\def\ninebig#1{{\hbox{$\textfont0=\tenrm\textfont2=\tensy
        \left#1\vbox to7.25pt{}\right.$}}}

%% TWELVE POINT FONT FAMILY --- not really small

\def\twelvepoint{\def\rm{\fam0\twelverm}
  \textfont0=\twelverm \scriptfont0=\eightrm \scriptscriptfont0=\sixrm
  \textfont1=\twelvei  \scriptfont1=\eighti  \scriptscriptfont1=\sixi
  \textfont2=\twelvesy  \scriptfont2=\eightsy  \scriptscriptfont2=\sixsy
  \textfont3=\tenex  \scriptfont3=\tenex  \scriptscriptfont3=\tenex
  \textfont\itfam=\twelveit  \def\it{\fam\itfam\twelveit}
  \textfont\slfam=\twelvesl  \def\sl{\fam\slfam\twelvesl}
  \textfont\ttfam=\twelvett  \def\tt{\fam\ttfam\twelvett}
  \textfont\bffam=\twelvebf  \scriptfont\bffam=\eightbf
    \scriptscriptfont\bffam=\sixbf  \def\bf{\fam\bffam\twelvebf}
  \tt  \ttglue=.5em plus.25em minus.15em
  \normalbaselineskip=11pt
  \setbox\strutbox=\hbox{\vrule height8pt depth3pt width0pt}
  \let\sc=\sevenrm  \let\big=\twelvebig \normalbaselines\rm}

\def\twelvebig#1{{\hbox{$\textfont0=\tenrm\textfont2=\tensy
        \left#1\vbox to7.25pt{}\right.$}}}
%%%%%%%%%%%%%%%%%%%%%%%%%%%%%%%%% end of smallfonts.tex %%%%%%%%%%%%%%%%%%%%%%%%%%%%%%%%
%%%%%%%%%%%%%%%%%%%%%%%%%%%%%%% including param.2 %%%%%%%%%%%%%%%%%%%%%%%%%%%%%%%
%param.2
\magnification=\magstep1
\def\firstpage{1}
\pageno=\firstpage
%%%%%%%%%%%%%%%%%%%%%%%%%%%%%%%%% end of param.2 %%%%%%%%%%%%%%%%%%%%%%%%%%%%%%%%
%%%%%%%%%%%%%%%%%%%%%%%%%%%%%%% including fonts.6 %%%%%%%%%%%%%%%%%%%%%%%%%%%%%%%
%fonts.6
\font\fiverm=cmr5
\font\sevenrm=cmr7
\font\sevenbf=cmbx7
\font\eightrm=cmr8
\font\eightbf=cmbx8
\font\ninerm=cmr9
\font\ninebf=cmbx9
\font\tenbf=cmbx10
\font\magtenbf=cmbx10 scaled\magstep1

%

%

%
%%%%%%%%%%%%%%%%%%%%%%%%%%%%%%%%% end of fonts.6 %%%%%%%%%%%%%%%%%%%%%%%%%%%%%%%%
%%%%%%%%%%%%%%%%%%%%%%%%%%%%%%% including smallfonts.tex %%%%%%%%%%%%%%%%%%%%%%%%%%%%%%%
%smallfonts.tex
%
\newskip\ttglue
\font\fiverm=cmr5
\font\fivei=cmmi5
\font\fivesy=cmsy5
\font\fivebf=cmbx5
\font\sixrm=cmr6
\font\sixi=cmmi6
\font\sixsy=cmsy6
\font\sixbf=cmbx6
\font\sevenrm=cmr7
\font\eightrm=cmr8
\font\eighti=cmmi8
\font\eightsy=cmsy8
\font\eightit=cmti8
\font\eightsl=cmsl8
\font\eighttt=cmtt8
\font\eightbf=cmbx8
\font\ninerm=cmr9
\font\ninei=cmmi9
\font\ninesy=cmsy9
\font\nineit=cmti9
\font\ninesl=cmsl9
\font\ninett=cmtt9
\font\ninebf=cmbx9
\font\twelverm=cmr12
\font\twelvei=cmmi12
\font\twelvesy=cmsy12
\font\twelveit=cmti12
\font\twelvesl=cmsl12
\font\twelvett=cmtt12
\font\twelvebf=cmbx12

%% EIGHT POINT FONT FAMILY

\def\eightpoint{\def\rm{\fam0\eightrm}
  \textfont0=\eightrm \scriptfont0=\sixrm \scriptscriptfont0=\fiverm
  \textfont1=\eighti  \scriptfont1=\sixi  \scriptscriptfont1=\fivei
  \textfont2=\eightsy  \scriptfont2=\sixsy  \scriptscriptfont2=\fivesy
  \textfont3=\tenex  \scriptfont3=\tenex  \scriptscriptfont3=\tenex
  \textfont\itfam=\eightit  \def\it{\fam\itfam\eightit}
  \textfont\slfam=\eightsl  \def\sl{\fam\slfam\eightsl}
  \textfont\ttfam=\eighttt  \def\tt{\fam\ttfam\eighttt}
  \textfont\bffam=\eightbf  \scriptfont\bffam=\sixbf
    \scriptscriptfont\bffam=\fivebf  \def\bf{\fam\bffam\eightbf}
  \tt  \ttglue=.5em plus.25em minus.15em
  \normalbaselineskip=9pt
  \setbox\strutbox=\hbox{\vrule height7pt depth2pt width0pt}
  \let\sc=\sixrm  \let\big=\eightbig \normalbaselines\rm}

\def\eightbig#1{{\hbox{$\textfont0=\ninerm\textfont2=\ninesy
        \left#1\vbox to6.5pt{}\right.$}}}

%% NINE POINT FONT FAMILY

\def\ninepoint{\def\rm{\fam0\ninerm}
  \textfont0=\ninerm \scriptfont0=\sixrm \scriptscriptfont0=\fiverm
  \textfont1=\ninei  \scriptfont1=\sixi  \scriptscriptfont1=\fivei
  \textfont2=\ninesy  \scriptfont2=\sixsy  \scriptscriptfont2=\fivesy
  \textfont3=\tenex  \scriptfont3=\tenex  \scriptscriptfont3=\tenex
  \textfont\itfam=\nineit  \def\it{\fam\itfam\nineit}
  \textfont\slfam=\ninesl  \def\sl{\fam\slfam\ninesl}
  \textfont\ttfam=\ninett  \def\tt{\fam\ttfam\ninett}
  \textfont\bffam=\ninebf  \scriptfont\bffam=\sixbf
    \scriptscriptfont\bffam=\fivebf  \def\bf{\fam\bffam\ninebf}
  \tt  \ttglue=.5em plus.25em minus.15em
  \normalbaselineskip=11pt
  \setbox\strutbox=\hbox{\vrule height8pt depth3pt width0pt}
  \let\sc=\sevenrm  \let\big=\ninebig \normalbaselines\rm}

\def\ninebig#1{{\hbox{$\textfont0=\tenrm\textfont2=\tensy
        \left#1\vbox to7.25pt{}\right.$}}}

%% TWELVE POINT FONT FAMILY --- not really small

\def\twelvepoint{\def\rm{\fam0\twelverm}
  \textfont0=\twelverm \scriptfont0=\eightrm \scriptscriptfont0=\sixrm
  \textfont1=\twelvei  \scriptfont1=\eighti  \scriptscriptfont1=\sixi
  \textfont2=\twelvesy  \scriptfont2=\eightsy  \scriptscriptfont2=\sixsy
  \textfont3=\tenex  \scriptfont3=\tenex  \scriptscriptfont3=\tenex
  \textfont\itfam=\twelveit  \def\it{\fam\itfam\twelveit}
  \textfont\slfam=\twelvesl  \def\sl{\fam\slfam\twelvesl}
  \textfont\ttfam=\twelvett  \def\tt{\fam\ttfam\twelvett}
  \textfont\bffam=\twelvebf  \scriptfont\bffam=\eightbf
    \scriptscriptfont\bffam=\sixbf  \def\bf{\fam\bffam\twelvebf}
  \tt  \ttglue=.5em plus.25em minus.15em
  \normalbaselineskip=11pt
  \setbox\strutbox=\hbox{\vrule height8pt depth3pt width0pt}
  \let\sc=\sevenrm  \let\big=\twelvebig \normalbaselines\rm}

\def\twelvebig#1{{\hbox{$\textfont0=\tenrm\textfont2=\tensy
        \left#1\vbox to7.25pt{}\right.$}}}
\catcode`\@=11
%
%  Include all definitions related to the fonts msam, msbm and eufm, so that
%  when this file is used by itself, the results with respect to those fonts
%  are equivalent to what they would have been using AMS-TeX.
%  Most symbols in fonts msam and msbm are defined using \newsymbol;
%  however, a few symbols that replace composites defined in plain must be
%  defined with \mathchardef.

\def\undefine#1{\let#1\undefined}
\def\newsymbol#1#2#3#4#5{\let\next@\relax
 \ifnum#2=\@ne\let\next@\msafam@\else
 \ifnum#2=\tw@\let\next@\msbfam@\fi\fi
 \mathchardef#1="#3\next@#4#5}
\def\mathhexbox@#1#2#3{\relax
 \ifmmode\mathpalette{}{\m@th\mathchar"#1#2#3}%
 \else\leavevmode\hbox{$\m@th\mathchar"#1#2#3$}\fi}
\def\hexnumber@#1{\ifcase#1 0\or 1\or 2\or 3\or 4\or 5\or 6\or 7\or 8\or
 9\or A\or B\or C\or D\or E\or F\fi}

\font\tenmsa=msam10
\font\sevenmsa=msam7
\font\fivemsa=msam5
\newfam\msafam
\textfont\msafam=\tenmsa
\scriptfont\msafam=\sevenmsa
\scriptscriptfont\msafam=\fivemsa
\edef\msafam@{\hexnumber@\msafam}
\mathchardef\dabar@"0\msafam@39
\def\dashrightarrow{\mathrel{\dabar@\dabar@\mathchar"0\msafam@4B}}
\def\dashleftarrow{\mathrel{\mathchar"0\msafam@4C\dabar@\dabar@}}

\def\ulcorner{\delimiter"4\msafam@70\msafam@70 }
\def\urcorner{\delimiter"5\msafam@71\msafam@71 }
\def\llcorner{\delimiter"4\msafam@78\msafam@78 }
\def\lrcorner{\delimiter"5\msafam@79\msafam@79 }
%    Note that there should not be a final space after the digits for a
%    \mathhexbox@.
\def\yen{{\mathhexbox@\msafam@55}}
\def\checkmark{{\mathhexbox@\msafam@58}}
\def\circledR{{\mathhexbox@\msafam@72}}
\def\maltese{{\mathhexbox@\msafam@7A}}

\font\tenmsb=msbm10
\font\sevenmsb=msbm7
\font\fivemsb=msbm5
\newfam\msbfam
\textfont\msbfam=\tenmsb
\scriptfont\msbfam=\sevenmsb
\scriptscriptfont\msbfam=\fivemsb
\edef\msbfam@{\hexnumber@\msbfam}
\def\Bbb#1{{\fam\msbfam\relax#1}}
\def\widehat#1{\setbox\z@\hbox{$\m@th#1$}%
 \ifdim\wd\z@>\tw@ em\mathaccent"0\msbfam@5B{#1}%
 \else\mathaccent"0362{#1}\fi}
\def\widetilde#1{\setbox\z@\hbox{$\m@th#1$}%
 \ifdim\wd\z@>\tw@ em\mathaccent"0\msbfam@5D{#1}%
 \else\mathaccent"0365{#1}\fi}
\font\teneufm=eufm10
\font\seveneufm=eufm7
\font\fiveeufm=eufm5
\newfam\eufmfam
\textfont\eufmfam=\teneufm
\scriptfont\eufmfam=\seveneufm
\scriptscriptfont\eufmfam=\fiveeufm

\catcode`\@=11
%%  Load amssym.def if necessary: If \newsymbol is undefined, do nothing
%%  and the following \input statement will be executed; otherwise
%%  change \input to a temporary no-op.
%#\ifx\undefined\newsymbol \else \begingroup\def\input#1 {\endgroup}\fi
%#\input amssym.def \relax
%%  Most symbols in fonts msam and msbm are defined using \newsymbol.  A few
%%  that are delimiters or otherwise require special treatment have already
%%  been defined as soon as the fonts were loaded.  Finally, a few symbols
%%  that replace composites defined in plain must be undefined first.
\newsymbol\boxdot 1200
\newsymbol\boxplus 1201
\newsymbol\boxtimes 1202
\newsymbol\square 1003
\newsymbol\blacksquare 1004
\newsymbol\centerdot 1205
\newsymbol\lozenge 1006
\newsymbol\blacklozenge 1007
\newsymbol\circlearrowright 1308
\newsymbol\circlearrowleft 1309
\undefine\rightleftharpoons
\newsymbol\rightleftharpoons 130A
\newsymbol\leftrightharpoons 130B
\newsymbol\boxminus 120C
\newsymbol\Vdash 130D
\newsymbol\Vvdash 130E
\newsymbol\vDash 130F
\newsymbol\twoheadrightarrow 1310
\newsymbol\twoheadleftarrow 1311
\newsymbol\leftleftarrows 1312
\newsymbol\rightrightarrows 1313
\newsymbol\upuparrows 1314
\newsymbol\downdownarrows 1315
\newsymbol\upharpoonright 1316
 
\newsymbol\downharpoonright 1317
\newsymbol\upharpoonleft 1318
\newsymbol\downharpoonleft 1319
\newsymbol\rightarrowtail 131A
\newsymbol\leftarrowtail 131B
\newsymbol\leftrightarrows 131C
\newsymbol\rightleftarrows 131D
\newsymbol\Lsh 131E
\newsymbol\Rsh 131F
\newsymbol\rightsquigarrow 1320
\newsymbol\leftrightsquigarrow 1321
\newsymbol\looparrowleft 1322
\newsymbol\looparrowright 1323
\newsymbol\circeq 1324
\newsymbol\succsim 1325
\newsymbol\gtrsim 1326
\newsymbol\gtrapprox 1327
\newsymbol\multimap 1328
\newsymbol\therefore 1329
\newsymbol\because 132A
\newsymbol\doteqdot 132B
 
\newsymbol\triangleq 132C
\newsymbol\precsim 132D
\newsymbol\lesssim 132E
\newsymbol\lessapprox 132F
\newsymbol\eqslantless 1330
\newsymbol\eqslantgtr 1331
\newsymbol\curlyeqprec 1332
\newsymbol\curlyeqsucc 1333
\newsymbol\preccurlyeq 1334
\newsymbol\leqq 1335
\newsymbol\leqslant 1336
\newsymbol\lessgtr 1337
\newsymbol\backprime 1038
\newsymbol\risingdotseq 133A
\newsymbol\fallingdotseq 133B
\newsymbol\succcurlyeq 133C
\newsymbol\geqq 133D
\newsymbol\geqslant 133E
\newsymbol\gtrless 133F
\newsymbol\sqsubset 1340
\newsymbol\sqsupset 1341
\newsymbol\vartriangleright 1342
\newsymbol\vartriangleleft 1343
\newsymbol\trianglerighteq 1344
\newsymbol\trianglelefteq 1345
\newsymbol\bigstar 1046
\newsymbol\between 1347
\newsymbol\blacktriangledown 1048
\newsymbol\blacktriangleright 1349
\newsymbol\blacktriangleleft 134A
\newsymbol\vartriangle 134D
\newsymbol\blacktriangle 104E
\newsymbol\triangledown 104F
\newsymbol\eqcirc 1350
\newsymbol\lesseqgtr 1351
\newsymbol\gtreqless 1352
\newsymbol\lesseqqgtr 1353
\newsymbol\gtreqqless 1354
\newsymbol\Rrightarrow 1356
\newsymbol\Lleftarrow 1357
\newsymbol\veebar 1259
\newsymbol\barwedge 125A
\newsymbol\doublebarwedge 125B
\undefine\angle
\newsymbol\angle 105C
\newsymbol\measuredangle 105D
\newsymbol\sphericalangle 105E
\newsymbol\varpropto 135F
\newsymbol\smallsmile 1360
\newsymbol\smallfrown 1361
\newsymbol\Subset 1362
\newsymbol\Supset 1363
\newsymbol\Cup 1264
 
\newsymbol\Cap 1265
 
\newsymbol\curlywedge 1266
\newsymbol\curlyvee 1267
\newsymbol\leftthreetimes 1268
\newsymbol\rightthreetimes 1269
\newsymbol\subseteqq 136A
\newsymbol\supseteqq 136B
\newsymbol\bumpeq 136C
\newsymbol\Bumpeq 136D
\newsymbol\lll 136E
 
\newsymbol\ggg 136F
 
\newsymbol\circledS 1073
\newsymbol\pitchfork 1374
\newsymbol\dotplus 1275
\newsymbol\backsim 1376
\newsymbol\backsimeq 1377
\newsymbol\complement 107B
\newsymbol\intercal 127C
\newsymbol\circledcirc 127D
\newsymbol\circledast 127E
\newsymbol\circleddash 127F
\newsymbol\lvertneqq 2300
\newsymbol\gvertneqq 2301
\newsymbol\nleq 2302
\newsymbol\ngeq 2303
\newsymbol\nless 2304
\newsymbol\ngtr 2305
\newsymbol\nprec 2306
\newsymbol\nsucc 2307
\newsymbol\lneqq 2308
\newsymbol\gneqq 2309
\newsymbol\nleqslant 230A
\newsymbol\ngeqslant 230B
\newsymbol\lneq 230C
\newsymbol\gneq 230D
\newsymbol\npreceq 230E
\newsymbol\nsucceq 230F
\newsymbol\precnsim 2310
\newsymbol\succnsim 2311
\newsymbol\lnsim 2312
\newsymbol\gnsim 2313
\newsymbol\nleqq 2314
\newsymbol\ngeqq 2315
\newsymbol\precneqq 2316
\newsymbol\succneqq 2317
\newsymbol\precnapprox 2318
\newsymbol\succnapprox 2319
\newsymbol\lnapprox 231A
\newsymbol\gnapprox 231B
\newsymbol\nsim 231C
\newsymbol\ncong 231D
\newsymbol\diagup 201E
\newsymbol\diagdown 201F
\newsymbol\varsubsetneq 2320
\newsymbol\varsupsetneq 2321
\newsymbol\nsubseteqq 2322
\newsymbol\nsupseteqq 2323
\newsymbol\subsetneqq 2324
\newsymbol\supsetneqq 2325
\newsymbol\varsubsetneqq 2326
\newsymbol\varsupsetneqq 2327
\newsymbol\subsetneq 2328
\newsymbol\supsetneq 2329
\newsymbol\nsubseteq 232A
\newsymbol\nsupseteq 232B
\newsymbol\nparallel 232C
\newsymbol\nmid 232D
\newsymbol\nshortmid 232E
\newsymbol\nshortparallel 232F
\newsymbol\nvdash 2330
\newsymbol\nVdash 2331
\newsymbol\nvDash 2332
\newsymbol\nVDash 2333
\newsymbol\ntrianglerighteq 2334
\newsymbol\ntrianglelefteq 2335
\newsymbol\ntriangleleft 2336
\newsymbol\ntriangleright 2337
\newsymbol\nleftarrow 2338
\newsymbol\nrightarrow 2339
\newsymbol\nLeftarrow 233A
\newsymbol\nRightarrow 233B
\newsymbol\nLeftrightarrow 233C
\newsymbol\nleftrightarrow 233D
\newsymbol\divideontimes 223E
\newsymbol\varnothing 203F
\newsymbol\nexists 2040
\newsymbol\Finv 2060
\newsymbol\Game 2061
\newsymbol\mho 2066
\newsymbol\eth 2067
\newsymbol\eqsim 2368
\newsymbol\beth 2069
\newsymbol\gimel 206A
\newsymbol\daleth 206B
\newsymbol\lessdot 236C
\newsymbol\gtrdot 236D
\newsymbol\ltimes 226E
\newsymbol\rtimes 226F
\newsymbol\shortmid 2370
\newsymbol\shortparallel 2371
\newsymbol\smallsetminus 2272
\newsymbol\thicksim 2373
\newsymbol\thickapprox 2374
\newsymbol\approxeq 2375
\newsymbol\succapprox 2376
\newsymbol\precapprox 2377
\newsymbol\curvearrowleft 2378
\newsymbol\curvearrowright 2379
\newsymbol\digamma 207A
\newsymbol\varkappa 207B
\newsymbol\Bbbk 207C
\newsymbol\hslash 207D
\undefine\hbar
\newsymbol\hbar 207E
\newsymbol\backepsilon 237F
%  Restore the catcode value for @ that was previously saved.
%#\catcode`\@=\csname pre amssym.tex at\endcsname

%\endinput
%%%%%%%%%%%%%%%%%%%%%%%%%%%%%%%%% end of symbols.1 %%%%%%%%%%%%%%%%%%%%%%%%%%%%%%%%
%%%%%%%%%%%%%%%%%%%%%%%%%%%%%%% including links.1 %%%%%%%%%%%%%%%%%%%%%%%%%%%%%%%
% links.1
% adapted from http://insti.physics.sunysb.edu/~siegel/tex.shtml
%
% postscript/pdf
\newcount\marknumber    \marknumber=1
\newcount\countdp \newcount\countwd \newcount\countht
%
% for ordinary tex
%
\ifx\pdfoutput\undefined
\def\rgboo#1{}
\def\postscript#1{\special{" #1}}               %% for dvips
\postscript{
        /bd {bind def} bind def
        /fsd {findfont exch scalefont def} bd
        /sms {setfont moveto show} bd
        /ms {moveto show} bd
        /pdfmark where          % printers ignore pdfmarks
        {pop} {userdict /pdfmark /cleartomark load put} ifelse
        [ /PageMode /UseOutlines                % bookmark window open
        /DOCVIEW pdfmark}
\def\bookmark#1#2{\postscript{          % #1=subheadings (if not 0)
        [ /Dest /MyDest\the\marknumber /View [ /XYZ null null null ] /DEST pdfmark
        [ /Title (#2) /Count #1 /Dest /MyDest\the\marknumber /OUT pdfmark}%
        \advance\marknumber by1}
\def\pdfclink#1#2#3{%
        \hskip-.25em\setbox0=\hbox{#2}%
                \countdp=\dp0 \countwd=\wd0 \countht=\ht0%
                \divide\countdp by65536 \divide\countwd by65536%
                        \divide\countht by65536%
                \advance\countdp by1 \advance\countwd by1%
                        \advance\countht by1%
                \def\linkdp{\the\countdp} \def\linkwd{\the\countwd}%
                        \def\linkht{\the\countht}%
        \postscript{
                [ /Rect [ -1.5 -\linkdp.0 0\linkwd.0 0\linkht.5 ]
                /Border [ 0 0 0 ]
                /Action << /Subtype /URI /URI (#3) >>
                /Subtype /Link
                /ANN pdfmark}{\rgb{#1}{#2}}}
%
% for pdftex
%
\else
\def\rgboo#1{\pdfliteral{#1 rg #1 RG}}
\pdfcatalog{/PageMode /UseOutlines}             % bookmark window open
\def\bookmark#1#2{
        \pdfdest num \marknumber xyz
        \pdfoutline goto num \marknumber count #1 {#2}
        \advance\marknumber by1}
\def\pdfklink#1#2{%
        \noindent\pdfstartlink user
                {/Subtype /Link
                /Border [ 0 0 0 ]
                /A << /S /URI /URI (#2) >>}{\rgb{1 0 0}{#1}}%
        \pdfendlink}
\fi

\def\rgbo#1#2{\rgboo{#1}#2\rgboo{0 0 0}}
\def\rgb#1#2{\mark{#1}\rgbo{#1}{#2}\mark{0 0 0}}
\def\pdfklink#1#2{\pdfclink{1 0 0}{#1}{#2}}
\def\pdflink#1{\pdfklink{#1}{#1}}
%
% examples:
% \bookmark{0}{look here}
% \pdfclink{0 0 1}{testlink}{http://www.google.com/}
% \pdfklink{testlink}{http://www.google.com/}
% \pdflink{http://www.google.com/}
%%%%%%%%%%%%%%%%%%%%%%%%%%%%%%%%% end of links.1 %%%%%%%%%%%%%%%%%%%%%%%%%%%%%%%%
%%%%%%%%%%%%%%%%%%%%%%%%%%%%%%% including titles.9 %%%%%%%%%%%%%%%%%%%%%%%%%%%%%%%
%titles.8
% requires fonts.5 or higher and smallfonts.tex
% uses links.* if included
% enumerates \demo consecutively (no section number)
%
\newcount\seccount  %% sections
\newcount\subcount  %% subsection
\newcount\clmcount  %% claim
\newcount\equcount  %% equation
\newcount\refcount  %% reference
\newcount\demcount  %% example
\newcount\execount  %% exercise
\newcount\procount  %% problem
\seccount=0
\equcount=1
\clmcount=1
\subcount=1
\refcount=1
\demcount=0
\execount=0
\procount=0
%
%% MISC STUFF
\def\proof{\medskip\noindent{\bf Proof.\ }}
\def\proofof(#1){\medskip\noindent{\bf Proof of \csname c#1\endcsname.\ }}
\def\qed{\hfill{\sevenbf QED}\par\medskip}
\def\references{\bigskip\noindent\hbox{\bf References}\medskip
                \ifx\pdflink\undefined\else\bookmark{0}{References}\fi}
\def\addref#1{\expandafter\xdef\csname r#1\endcsname{\number\refcount}
    \global\advance\refcount by 1}

\def\nextremark #1\par{\item{$\circ$} #1}
\def\firstremark #1\par{\bigskip\noindent{\bf Remarks.}
     \smallskip\nextremark #1\par}
\def\abstract#1\par{{\baselineskip=10pt
    \eightpoint\narrower\noindent{\eightbf Abstract.} #1\par}}
%
%% EQUATION
\def\equtag#1{\expandafter\xdef\csname e#1\endcsname{(\number\seccount.\number\equcount)}
              \global\advance\equcount by 1}
\def\equation(#1){\equtag{#1}\eqno\csname e#1\endcsname}
\def\equ(#1){\hskip-0.03em\csname e#1\endcsname}
%
%% CLAIMS (theorems etc)
\def\clmtag#1#2{\expandafter\xdef\csname cn#2\endcsname{\number\seccount.\number\clmcount}
                \expandafter\xdef\csname c#2\endcsname{#1~\number\seccount.\number\clmcount}
                \global\advance\clmcount by 1}
\def\claim #1(#2) #3\par{\clmtag{#1}{#2}
    \vskip.1in\medbreak\noindent
    {\bf \csname c#2\endcsname .\ }{\sl #3}\par
    \ifdim\lastskip<\medskipamount
    \removelastskip\penalty55\medskip\fi}
\def\clm(#1){\csname c#1\endcsname}
\def\clmno(#1){\csname cn#1\endcsname}
%
%% SECTION
\def\sectag#1{\global\advance\seccount by 1
              \expandafter\xdef\csname sectionname\endcsname{\number\seccount. #1}
              \equcount=1 \clmcount=1 \subcount=1 \execount=0 \procount=0}
\def\section#1\par{\vskip0pt plus.1\vsize\penalty-40
    \vskip0pt plus -.1\vsize\bigskip\bigskip
    \sectag{#1}
    \message{\sectionname}\leftline{\magtenbf\sectionname}
    \nobreak\smallskip\noindent
    \ifx\pdflink\undefined
    \else
      \bookmark{0}{\sectionname}
    \fi}
%
%% SUBSECTION
\def\subtag#1{\expandafter\xdef\csname subsectionname\endcsname{\number\seccount.\number\subcount. #1}
              \global\advance\subcount by 1}
\def\subsection#1\par{\vskip0pt plus.05\vsize\penalty-20
    \vskip0pt plus -.05\vsize\medskip\medskip
    \subtag{#1}
    \message{\subsectionname}\leftline{\tenbf\subsectionname}
    \nobreak\smallskip\noindent
    \ifx\pdflink\undefined
    \else
      \bookmark{0}{.... \subsectionname}  %% can get a bit cluttered
    \fi}
%
%% DEMO (examples etc)
\def\demtag#1#2{\global\advance\demcount by 1
              \expandafter\xdef\csname de#2\endcsname{#1~\number\demcount}}
\def\demo #1(#2) #3\par{
  \demtag{#1}{#2}
  \vskip.1in\medbreak\noindent
  {\bf #1 \number\demcount.\enspace}
  {\rm #3}\par
  \ifdim\lastskip<\medskipamount
  \removelastskip\penalty55\medskip\fi}
\def\dem(#1){\csname de#1\endcsname}
%
%% EXERCISE
\def\exetag#1{\global\advance\execount by 1
              \expandafter\xdef\csname ex#1\endcsname{Exercise~\number\seccount.\number\execount}}
\def\exercise(#1) #2\par{
  \exetag{#1}
  \vskip.1in\medbreak\noindent
  {\bf Exercise \number\execount.}
  {\rm #2}\par
  \ifdim\lastskip<\medskipamount
  \removelastskip\penalty55\medskip\fi}
\def\exe(#1){\csname ex#1\endcsname}
%
%% PROBLEM
\def\protag#1{\global\advance\procount by 1
              \expandafter\xdef\csname pr#1\endcsname{\number\seccount.\number\procount}}
\def\problem(#1) #2\par{
  \ifnum\procount=0
    \parskip=6pt
    \vbox{\bigskip\centerline{\bf Problems \number\seccount}\nobreak\medskip}
  \fi
  \protag{#1}
  \item{\number\procount.} #2}
\def\pro(#1){Problem \csname pr#1\endcsname}
%
%%%%%%%%%%%%%%%%%%%%%%%%%%%%%%%%% end of titles.9 %%%%%%%%%%%%%%%%%%%%%%%%%%%%%%%%
%%%%%%%%%%%%%%%%%%%%%%%%%%%%%%% including macros.21 %%%%%%%%%%%%%%%%%%%%%%%%%%%%%%%
%macros.20
%
% requires fonts.5 or later
% also defines mathds (double strike) family
%
\def\rightheadline{\hfil}
\def\leftheadline{\sevenrm\hfil HANS KOCH\hfil}
\headline={\ifnum\pageno=\firstpage\hfil\else
\ifodd\pageno{{\fiverm\rightheadline}\number\pageno}
\else{\number\pageno\fiverm\leftheadline}\fi\fi}
\footline={\ifnum\pageno=\firstpage\hss\tenrm\folio\hss\else\hss\fi}
\let\ov=\overline
\let\cl=\centerline

\let\eps=\varepsilon
\let\sss=\scriptscriptstyle

\def\AA{{\cal A}}
\def\BB{{\cal B}}

\def\DD{{\cal D}}
\def\EE{{\cal E}}

\def\HH{{\cal H}}

\def\LL{{\cal L}}

\def\NN{{\cal N}}

\def\PP{{\cal P}}

\def\RR{{\cal R}}
\def\SS{{\cal S}}

\def\XX{{\cal X}}
\def\YY{{\cal Y}}

\def\ssL{{\sss L}}

\def\ssR{{\sss R}}

\def\rmC{\mathop{\rm C}\nolimits}
\def\rmL{\mathop{\rm L}\nolimits}
\def\id{\mathop{\rm I}\nolimits}

\def\Im{\mathop{\rm Im}\nolimits}
%
%%%%%%%%%%%%%%
\newfam\dsfam
\def\mathds #1{{\fam\dsfam\tends #1}}

\font\tends=dsrom10
\font\eightds=dsrom8
\textfont\dsfam=\tends
\scriptfont\dsfam=\eightds
%%%%%%%%%%%%%%
%

\def\integer{{\mathds Z}}

\def\real{{\mathds R}}

\def\proj{{\Bbb P}}

\def\bdot{\hbox{\bf .}}

\def\defeq{\mathrel{\mathop=^{\sss\rm def}}}
\def\half{{1\over 2}}

\def\thalf{{\textstyle\half}}

%

%

%

%

%
% from TeX book: used for commutative diagram
% in math mode, before using matrix, do
% \def\normalbaselines{\baselineskip20pt\lineskip3pt\lineskiplimit3pt}

%%%%%%%%%%%%%%%%%%%%%%%%%%%%%%%%% end of macros.21 %%%%%%%%%%%%%%%%%%%%%%%%%%%%%%%%
%%%%%%%%%%%%%%%%%%%%%%%%%%%%%%% including mygraphicx.tex %%%%%%%%%%%%%%%%%%%%%%%%%%%%%%%
%% modification of graphicx.tex by Nathan Goldschmidt
\input miniltx

\ifx\pdfoutput\undefined
  \def\Gin@driver{dvips.def}  % we are not running PDFTeX
\else
  \def\Gin@driver{pdftex.def} % we are running PDFTeX
\fi

\input graphicx.sty
\resetatcatcode
%%%%%%%%%%%%%%%%%%%%%%%%%%%%%%%%% end of mygraphicx.tex %%%%%%%%%%%%%%%%%%%%%%%%%%%%%%%%
%%%%%%%%%%%%%%%%%%%%%%%%%%%%%%% including opmac.1 %%%%%%%%%%%%%%%%%%%%%%%%%%%%%%%
%% table macros from opmac.tex
%% Petr Olsak, 2012 -- 2016
%% http://petr.olsak.net/opmac.html

\newcount\tmpnum % auxiliary count
\newdimen\tmpdim % auxiliary dimen
\def\opwarning#1{\immediate\write16{l.\the\inputlineno\space OPmac WARNING: #1.}}
\long\def\addto#1#2{\expandafter\def\expandafter#1\expandafter{#1#2}}
\long\def\isinlist#1#2#3{\begingroup \long\def\tmp##1#2##2\end{\def\tmp{##2}%
   \ifx\tmp\empty \endgroup \csname iffalse\expandafter\endcsname \else
                  \endgroup \csname iftrue\expandafter\endcsname \fi}% end of \def\tmp
   \expandafter\tmp#1\endlistsep#2\end
}
\def\tabstrut{\strut}     % strut in the \table
\def\tabiteml{\enspace}   % left material before each \table item
\def\tabitemr{\enspace}   % right material after each \table item
\def\vvkern{1pt}          % space between vertical lines
\def\hhkern{1pt}          % space between horizontal lines

%%%%%%%%%%%%%% \table -- sec. 3.19 in opmac-d.pdf

\newtoks\tabdata
\def\tabstrutA{\tabstrut}
\newcount\colnum
\def\ddlinedata{}
\def\vvleft{}

\def\table{\vbox\bgroup \catcode`\|=12 \tableA}
\def\tableA#1#2{\offinterlineskip \colnum=0 \def\tmpa{}\tabdata={}\scantabdata#1\relax
   \halign\expandafter{\the\tabdata\cr#2\crcr}\egroup}

\def\scantabdata#1{\let\next=\scantabdata
   \ifx\relax#1\let\next=\relax
   \else\ifx|#1\addtabvrule
      \else\isinlist{123456789}#1\iftrue \def\next{\scantabdataC#1}%
          \else \expandafter\ifx\csname tabdeclare#1\endcsname \relax
                \expandafter\ifx\csname paramtabdeclare#1\endcsname \relax
                   \opwarning{tab-declarator "#1" unknown, ignored}%
                \else \def\next{\expandafter \scantabdataB \csname paramtabdeclare#1\endcsname}\fi
             \else \def\next{\expandafter\scantabdataA \csname tabdeclare#1\endcsname}%
   \fi\fi\fi\fi \next
}
\def\scantabdataA#1{\addtabitem \expandafter\addtabdata\expandafter{#1\tabstrutA}\scantabdata}
\def\scantabdataB#1#2{\addtabitem\expandafter\addtabdata\expandafter{#1{#2}\tabstrutA}\scantabdata}
\def\scantabdataC {\def\tmpb{}\afterassignment\scantabdataD \tmpnum=}
\def\scantabdataD#1{\loop \ifnum\tmpnum>0 \advance\tmpnum by-1 \addto\tmpb{#1}\repeat
   \expandafter\scantabdata\tmpb
}

\def\unsskip{\ifdim\lastskip>0pt \unskip\fi}
\def\addtabitem{\ifnum\colnum>0 \addtabdata{&}\addto\ddlinedata{&\dditem}\fi
    \advance\colnum by1 \let\tmpa=\relax}
\def\addtabdata#1{\tabdata\expandafter{\the\tabdata#1}}
\def\addtabvrule{%
    \ifx\tmpa\vrule \addtabdata{\kern\vvkern}%
       \ifnum\colnum=0 \addto\vvleft{\vvitem}\else\addto\ddlinedata{\vvitem}\fi
    \else \ifnum\colnum=0 \addto\vvleft{\vvitemA}\else\addto\ddlinedata{\vvitemA}\fi\fi
    \let\tmpa=\vrule \addtabdata{\vrule}}

\def\crl{\crcr\noalign{\hrule}}
\def\crll{\crcr\noalign{\hrule\kern\hhkern\hrule}}

\def\crli{\crcr \omit
   \gdef\dditem{\omit\tablinefil}\gdef\vvitem{\kern\vvkern\vrule}\gdef\vvitemA{\vrule}%
   \vvleft\tablinefil\ddlinedata\crcr}
\def\crlli{\crli\noalign{\kern\hhkern}\crli}
\def\tablinefil{\leaders\hrule\hfil}

\def\crlp#1{\crcr \noalign{\kern-\drulewidth}%
   \omit \xdef\crlplist{#1}\xdef\crlplist{,\expandafter}\expandafter\crlpA\crlplist,\end,%
   \global\tmpnum=0 \gdef\dditem{\omit\crlpD}%
   \gdef\vvitem{\kern\vvkern\kern\drulewidth}\gdef\vvitemA{\kern\drulewidth}%
   \vvleft\crlpD\ddlinedata \global\tmpnum=0 \crcr}
\def\crlpA#1,{\ifx\end#1\else \crlpB#1-\end,\expandafter\crlpA\fi}
\def\crlpB#1#2-#3,{\ifx\end#3\xdef\crlplist{\crlplist#1#2,}\else\crlpC#1#2-#3,\fi}
\def\crlpC#1-#2-#3,{\tmpnum=#1\relax
   \loop \xdef\crlplist{\crlplist\the\tmpnum,}\ifnum\tmpnum<#2\advance\tmpnum by1 \repeat}
\def\crlpD{\global\advance\tmpnum by1
   \edef\tmpa{\noexpand\isinlist\noexpand\crlplist{,\the\tmpnum,}}%
   \tmpa\iftrue \kern-\drulewidth \tablinefil \kern-\drulewidth\else\hfil \fi}

\def\tskip{\afterassignment\tskipA \tmpdim}
\def\tskipA{\gdef\dditem{}\gdef\vvitem{}\gdef\vvitemA{}\gdef\tabstrutA{}%
    \vbox to\tmpdim{}\ddlinedata \crcr \noalign{\gdef\tabstrutA{\tabstrut}}}

\def\mspan{\omit \tabdata={\tabstrut}\let\tmpa=\relax \afterassignment\mspanA \mscount=}
\def\mspanA[#1]{\loop \ifnum\mscount>1 \csname span\endcsname \omit \advance\mscount by-1 \repeat
   \mspanB#1\relax}
\def\mspanB#1{\ifx\relax#1\def\tmpa{\def\tmpa####1}%
   \expandafter\tmpa\expandafter{\the\tabdata\ignorespaces}\expandafter\tmpa\else
   \ifx |#1\ifx\tmpa\vrule\addtabdata{\kern\vvkern}\fi \addtabdata{\vrule}\let\tmpa=\vrule
   \else \let\tmpa=\relax
      \ifx c#1\addtabdata{\tabiteml\hfil\ignorespaces##1\unsskip\hfil\tabitemr}\fi
      \ifx l#1\addtabdata{\tabiteml\ignorespaces##1\unsskip\hfil\tabitemr}\fi
      \ifx r#1\addtabdata{\tabiteml\hfil\ignorespaces##1\unsskip\tabitemr}\fi
   \fi \expandafter\mspanB \fi}

\newdimen\drulewidth  \drulewidth=0.4pt
\let\orihrule=\hrule  \let\orivrule=\vrule
\def\rulewidth{\afterassignment\rulewidthA \drulewidth}
\def\rulewidthA{\edef\hrule{\orihrule height\the\drulewidth}%
                \edef\vrule{\orivrule width\the\drulewidth}}

\long\def\frame#1{%
   \hbox{\vrule\vtop{\vbox{\hrule\kern\vvkern
      \hbox{\kern\hhkern\relax#1\kern\hhkern}%
   }\kern\vvkern\hrule}\vrule}}
%%%%%%%%%%%%%%%%%%%%%%%%%%%%%%%%% end of opmac.1 %%%%%%%%%%%%%%%%%%%%%%%%%%%%%%%%
%%%%%%%%%%%%%%%%%%%%%%%%%%%%%%% including boldmath.tex %%%%%%%%%%%%%%%%%%%%%%%%%%%%%%%
%% boldmath

%=================================================================
% necessary fonts
%-----------------------------------------------------------------
   \font\Fivebf  =cmbx10  scaled 500 % five point bold
   \font\Sevenbf =cmbx10  scaled 700 % seven point bold
   \font\Tenbf   =cmbx10             % ten point bold
   \font\Fivemb  =cmmib10 scaled 500 % five point math bold
   \font\Sevenmb =cmmib10 scaled 700 % seven point math bold
   \font\Tenmb   =cmmib10            % ten point math bold
%=================================================================
% Math in Bold
% example $\boldmath{..}$   { math characters will be bold }
%-----------------------------------------------------------------
\def\boldmath{\textfont0=\Tenbf           \scriptfont0=\Sevenbf
              \scriptscriptfont0=\Fivebf  \textfont1=\Tenmb
              \scriptfont1=\Sevenmb       \scriptscriptfont1=\Fivemb}
%=================================================================

%% twelveboldmath

%=================================================================
% necessary fonts
%-----------------------------------------------------------------
%  \font\Sevenbf   =cmbx10  scaled 700

%  \font\Sevenmb   =cmmib10 scaled 700

%=================================================================
% Math in Bold
% example $\twelveboldmath{...}$   { math characters will be bold }
%-----------------------------------------------------------------

%=================================================================

%%%%%%%%%%%%%%%%%%%%%%%%%%%%%%%%% end of boldmath.tex %%%%%%%%%%%%%%%%%%%%%%%%%%%%%%%%
%\input param.2
%\input fonts.6
%\input smallfonts.tex
%\input symbols.1
%\input links.1
%\input titles.9
%\input macros.21
%\input boldmath.tex
%\input mygraphicx.tex
%\input opmac.1
%%%%%%%%%%%%%%%%%%%%%%% begin drawing macros
\input texdraw
\newcount\xval
\def\PPP{\advance\xval by 3\lvec({\number\xval}  3.0)\fcir f:0 r:0.6 }
\def\ppp{\advance\xval by 3\lvec({\number\xval}  1.5)\fcir f:0 r:0.6 }
\def\nnn{\advance\xval by 3\lvec({\number\xval} -1.5)\fcir f:0 r:0.6 }
\def\NNN{\advance\xval by 3\lvec({\number\xval} -3.0)\fcir f:0 r:0.6 }
\def\zzz{\advance\xval by 3\lvec({\number\xval}  0.0)}
\def\ZZZ{\advance\xval by 5\lvec({\number\xval}  0.0)}
\def\DRAW#1{\hbox{
\btexdraw\drawdim {pt}
\linewd 0.1 \xval=0
\ZZZ #1 \zzz\ZZZ\etexdraw}}
\def\SHOW#1{\hbox{
\btexdraw\drawdim {pt}
\linewd 0.1 \xval=0
#1 \zzz\etexdraw}}
%%%%%%%%%%%%%%%%%%%%% end drawing macros
\font\tenamsb=msbm10 \font\sevenamsb=msbm7 \font\fiveamsb=msbm5
\newfam\bbfam
\textfont\bbfam=\tenamsb
\scriptfont\bbfam=\sevenamsb
\scriptscriptfont\bbfam=\fiveamsb
\def\bbb{\fam\bbfam}
\def\oldinteger{{\bbb Z}}
\def\oldnatural{{\bbb N}}
\def\oldreal{{\bbb R}}

\let\Scale\kappa
\def\id{{\rm I}}
\def\mod{{\rm mod~}}
\def\col{{:\hskip3pt}}
\def\scol{{;\hskip3pt}}
\def\dotdot{\hskip1pt{..}\hskip1pt}
\def\ssL{{\sss L}}
\def\ssR{{\sss R}}
\def\even{0}
\def\odd{1}
\def\tinyskip{\hskip.7pt}
\def\Brr{\BB_{\rho,r}}
\def\Brrt{\Brr^{\,\tau}}
\def\bdot{\hbox{\bf .}}

\def\sfrac#1#2{\hbox{\raise2.2pt\hbox{$\scriptstyle#1$}\hskip-1.2pt
   {$\scriptstyle/$}\hskip-0.9pt\lower2.2pt\hbox{$\scriptstyle#2$}\hskip1.0pt}}
\def\shalf{\sfrac{1}{2}}
\addref{CBH}
\addref{Adams}
\addref{FrWa}
\addref{SmWi}
\addref{FPi}
\addref{IosK}
\addref{Ios}
\addref{FPii}
\addref{BuDe}
\addref{FPiii}
\addref{Zg}
\addref{APankov}
\addref{JaSii}
\addref{GGKKMO}
\addref{KaAl}
\addref{MiNa}
\addref{HeRa}
\addref{AKi}
\addref{PaRo}
\addref{HoWa}
\addref{AKii}
\addref{CZg}
\addref{SZg}
\addref{KLM}
\addref{KiLe}
\addref{Matsue}
\addref{AKiv}
\addref{AKv}
\addref{WZg}
\addref{Files}
\addref{Ada}
\addref{Gnat}
\addref{IEEE}
\addref{MPFR}
\def\leftheadline{\sixrm\hfil ARIOLI \& KOCH\hfil}
\def\rightheadline{\sevenrm\hfil traveling waves for FPU\hfil}
%
%%%%%%%%%%%%%%%%%%%%%%%%%%%%%%%%%%%%%%%%%%%%%%%%%%%%%%%%%%%%%%%%
\cl{{\magtenbf Traveling wave solutions for the FPU chain:}}
\cl{{\magtenbf a constructive approach}}
\bigskip

\cl{
Gianni Arioli
\footnote{$^1$}
{\eightpoint\hskip-2.9em
Department of Mathematics, Politecnico di Milano,
Piazza Leonardo da Vinci 32, 20133 Milano.
}
$^{\!\!\!,\!\!}$
\footnote{$^2$}
{\eightpoint\hskip-2.6em
Supported in part by the PRIN project
``Equazioni alle derivate parziali di tipo ellittico e parabolico: aspetti geometrici, disuguaglianze collegate, e applicazioni''.
}
and Hans Koch
\footnote{$^3$}
{\eightpoint\hskip-2.7em
Department of Mathematics, The University of Texas at Austin,
Austin, TX 78712.}
}

\bigskip
\abstract
Traveling waves for the FPU chain
are constructed by solving the associated equation
for the spatial profile $u$ of the wave.
We consider solutions whose derivatives $u'$
need not be small, may change sign several times,
but decrease at least exponentially.
Our method of proof is computer-assisted.
Unlike other methods,
it does not require that the FPU potential
has an attractive (positive) quadratic term.
But we currently need to restrict the size of that term.
In particular, our solutions in the attractive case
are all supersonic.

%%%%%%%%%%%%%%%%%%%%%
%%%%%%%%%%%%%%%%%%%%%
\section Introduction
%%%%%%%%%%%%%%%%%%%%%
%%%%%%%%%%%%%%%%%%%%%

We consider a chain of interacting particles
described by the equation
$$
\ddot q_j=\phi'(q_{j+1}-q_j)-\phi'(q_j-q_{j-1})\,,
\qquad j\in\integer\,,
\equation(FPU)
$$
where $\phi$ is a polynomial of degree at least $3$.
The choice
$$
\phi(v)={1\over 2}\phi_2 v^2+{1\over m+1}\phi_{m+1} v^{m+1}
\equation(phiv)
$$
corresponds to the FPU model:
the $\alpha$-model if $m=2$, or the $\beta$-model if $m=3$.

Our goal is to find traveling waves,
meaning solutions of the form $q_j(t)=u(j-t/\tau)$.
Substituting this ansatz into \equ(FPU) yields the equation
$$
u''(x)=\tau^2\phi'\bigl(u(x+1)-u(x))-\tau^2\phi'\bigl(u(x)-u(x-1)\bigr)\,,
\qquad x\in\real\,.
\equation(travelOne)
$$
We focus on solutions $u$ that approach limit values at $\pm\infty$.

The first result on the existence of traveling waves
for infinite chains of FPU type was obtained in [\rFrWa],
where solutions with prescribed energy are found as
constrained minima of a suitable functional.
This result has been extended to solutions with
prescribed speed in [\rSmWi], using the mountain pass theorem.
A survey of variational results, covering both breathers and traveling waves,
and including an extensive bibliography, is given in [\rAPankov].
Some recent results based on variational techniques
can be found e.g.~in [\rHeRa,\rPaRo].

A perturbative approach, based on center manifold theory,
has been developed in [\rIosK,\rIos,\rJaSii].
It yields small-amplitude solutions in cases where $\tau^2\phi_2\simeq 1$.
A different type of perturbative approach exploits the fact that,
for slowly varying functions, the equation \equ(FPU) is close to integrable.
This property has been used e.g.~in [\rFPi,\rFPii,\rFPiii,\rHoWa]
to construct and analyze solitary waves on FPU lattices.

These results all concern potentials
that are attractive for small displacements, meaning $\phi_2>0$.
In the approach considered here,
the value of $\phi_2$ is allowed to be negative or zero as well,
and we find several traveling wave solutions in this case.
A similar approach might apply to other equations of this type.
Advance-delay equations appear e.g.~in models from
biology [\rCBH], economics [\rKaAl], and electrodynamics [\rBuDe].

The method developed in this paper is constructive
but not limited to small or near-integrable solutions.
Starting with an approximate numerical solution,
we use computer-assisted methods to prove that
there exists a true solution nearby.
Our current technique requires that $|\tau^2\phi_2|<1$,
but we expect that this condition can be relaxed in future work.

Apart from this restriction,
our method applies in principle to any polynomial potential $\phi$.
But for simplicity, we restrict to the FPU model \equ(phiv).
Using as parameters $\mu=\tau^2\phi_2$ and $\nu=\tau^2\phi_{m+1}$,
we have
$$
\tau^2\phi'(v)=\mu v+\nu v^m\,.
\equation(muvnuvm)
$$
Without loss of generality, we may assume that $|\nu|=1$.
In addition, we restrict to solutions $u$ that are either even or odd.
In other words, $u(-x)=(-1)^\sigma u(x)$ for all $x$,
with $\sigma\in\{0,1\}$.
The number $\sigma$ will be referred to as the parity of $u$.

For any function $u$ on $\real$ define
$$
(\DD u)(x)=u\bigl(x+\thalf\bigr)-u\bigl(x-\thalf\bigr)\,.
\equation(DDuDef)
$$
Notice that $\DD u$ has parity $1-\sigma$, if $u$ has parity $\sigma$.

\claim Theorem(main)
Let $\nu=1$.
Consider a fixed but arbitrary row in Table 1
and the data given in that row.
Then for every value of $\mu$
in some open neighborhood of $\bar\mu$,
the equation \equ(travelOne),
with the potential given by \equ(muvnuvm),
has a real analytic solution $u$ with parity $\sigma$.
The function $v=\DD u$ satisfies a bound
$|v(x)|\le Cr^{-\kappa|x|}$ for some constant $C>0$.
Its $\sup$-norm is given in column 5,
where $p$ is some positive real number.
This function $v$ has $E$ local extrema with values $|v(t)|>\sfrac{1}{64}$.
(The existence of others is not excluded.)
The diagram in column 7 specifies the sequence and nature of these extrema,
as described below.

Each diagram in Table 1 represent the graph
of the function $v=\DD u$ associated with the solution $u$,
where the endpoints correspond to $x=\pm\infty$.
The vertices ``$\!\raise 2pt\SHOW{\PPP}\tinyskip$'',
``$\!\SHOW{\ppp}\tinyskip$'',
``$\!\raise 2pt\SHOW{\nnn}\tinyskip$'',
and ``$\!\raise-1pt\SHOW{\NNN}\tinyskip$'' represent
positive local maxima, negative local maxima, positive local minima,
and negative local minima, respectively;
and these extrema appear in the indicated order.
More detailed (but purely numerical) graphs are shown in Figures 1--8.

The parameters $r$, $k$, and $\Scale$ that are listed in Table 1
are used during our construction of the solution,
as will be explained later.
Concerning our choice for the values of $\mu$,
we note that similar solutions exist for many
(if not all) other values in the interval $(-1,0)$ or $(-1,1)$,
depending on whether $\mu$ is negative or nonnegative, respectively.
Further details will be given at the end of this section.

\bigskip
If $u$ is a solution of the equation \equ(travelOne),
then the function $v=\DD u$ satisfies
$$
v''=\tau^2\DD^2\phi'(v)\,.
\equation(travelTwo)
$$
We prove \clm(main) by first solving this equation
and verifying that the solution $v$ has the indicated properties.
This involves estimates that are verified with the aid of a computer.
Then we define two function $u_\ssL$ and $u_\ssR$ by setting
$$
u_\ssL(x)=\sum_{j=0}^\infty v\bigl(x-j-\thalf\bigr)\,,\qquad
u_\ssR(x)=-\sum_{j=0}^\infty v\bigl(x+j+\thalf\bigr)\,,
\equation(uLuR)
$$
for all $x\in\real$.
It is not hard to see that both $u=u_\ssL$ and $u=u_\ssR$
satisfy the equation \equ(travelOne),
and that the difference $u_\ssR-u_\ssL$ is constant.
By construction, both satisfy $\DD u=v$.
The solution with the proper parity is $u=\shalf u_\ssL+\shalf u_\ssR$.

%%%%%%%%%%%%%%%%%%%%%%%%%%%%%%%%%%%%%%%%%%%%%%%%%%%
%%%%%%%%%%%%%%%%%%%%%%%%%%%%%%%%%%%%%%%%%%%%%%%%%%% Table 1
\newif\ifprefix
%%\prefixtrue %% use prefix
\prefixfalse %% use number
\def\EonePZ{\ifprefix E1PZ\else 1\fi}
\def\EtwoNZ{\ifprefix E2NZ\else 2\fi}
\def\EonezeroA{\ifprefix E10A\else 3\fi}
\def\EonePA{\ifprefix E1PA\else 4\fi}
\def\EonePB{\ifprefix E1PB\else 5\fi}
\def\EoneNA{\ifprefix E1NA\else 6\fi}
\def\EoneNC{\ifprefix E1NC\else 7\fi}
\def\EtwoNA{\ifprefix E2NA\else 8\fi}
\def\EtwoNB{\ifprefix E2NB\else 9\fi}
\def\EtwoNC{\ifprefix E2NC\else 10\fi}
\def\EthreeNA{\ifprefix E3NA\else 11\fi}
\def\EthreeNB{\ifprefix E3NB\else 12\fi}
\def\EfourNA{\ifprefix E4NA\else 13\fi}
\def\EfourNB{\ifprefix E4NB\else 14\fi}
\def\EfiveNA{\ifprefix E5NA\else 15\fi}
\def\EfiveNB{\ifprefix E5NB\else 16\fi}
\def\EfiveNC{\ifprefix E5NC\else 17\fi}
\def\EfiveND{\ifprefix E5ND\else 18\fi}
\def\OtwoNA{\ifprefix O2NA\else 19\fi}
\def\OfourNA{\ifprefix O4NA\else 20\fi}
\def\OfourNB{\ifprefix O4NB\else 21\fi}
\def\OsixNA{\ifprefix O6NA\else 22\fi}
%%%%%%%%%%%%%%%%%%%%%%%%%%%%%%%%%%%%%%%%%%%%%%%%%%%
\bigskip\bigskip
%% l=left, c=center, r=right
\hfil\table{|c||c|c|c||c|c|c||c|c|c|}{\crl\tskip.1em
label & $m$ &$\bar\mu$&$\sigma$&$\|v\|_\infty$& $E$& diagram &$r$ &$k$&$\Scale$\crli\tskip.1em
\EonePZ & $2$ &  $1/4$& $1$ &$1.0+p$&$1$&\DRAW{\PPP}&$4$&$1$&$1$ \cr
\EtwoNZ & $2$ & $-1/4$& $1$ &$1.7+p$&$3$&\DRAW{\PPP\ppp\PPP}&$4$&$1$&$1$ \cr
\EonezeroA & $3$ &    $0$& $1$ &$1.3+p$&$1$&\DRAW{\PPP}&$4$&$1$&$1$ \cr
\EonePA & $3$ &  $1/2$& $1$ &$0.9+p$&$1$&\DRAW{\PPP}&$2$&$3$&$1$ \cr
\EonePB & $3$ &  $3/4$& $1$ &$0.6+p$&$1$&\DRAW{\PPP}&$3/2$&$8$&$1$ \cr
\EoneNA & $3$ & $-1/4$& $1$ &$1.4+p$&$1$&\DRAW{\PPP}&$4$&$2$&$1$ \cr
\EoneNC & $3$ & $-1/2$& $1$ &$1.5+p$&$1$&\DRAW{\PPP}&$2$&$2$&$2$ \cr
\EtwoNA & $3$ & $-1/4$& $1$ &$1.4+p$&$3$&\DRAW{\PPP\ppp\PPP}&$4$&$2$&$1$ \cr
\EtwoNB & $3$ &$-1/256$&$1$ &$1.3+p$&$3$&\DRAW{\PPP\ppp\PPP}&$3/2$&$1$&$1$ \cr
\EtwoNC & $3$ & $-3/4$& $1$ &$1.6+p$&$5$&\DRAW{\NNN\PPP\ppp\PPP\NNN}&$3/2$&$8$&$1$ \cr
\EthreeNA & $3$ & $-1/4$& $1$ &$1.4+p$&$5$&\DRAW{\PPP\ppp\PPP\ppp\PPP}&$4$&$2$&$1$ \cr
\EthreeNB & $3$ & $-1/4$& $1$ &$1.4+p$&$3$&\DRAW{\NNN\PPP\NNN}&$4$&$2$&$1$ \cr
\EfourNA & $3$ & $-1/4$& $1$ &$1.4+p$&$5$&\DRAW{\NNN\PPP\ppp\PPP\NNN}&$2$&$2$&$1$ \cr
\EfourNB & $3$ & $-1/4$& $1$ &$1.4+p$&$7$&\DRAW{\PPP\ppp\PPP\ppp\PPP\ppp\PPP}&$3$&$2$&$1$ \cr
\EfiveNA & $3$ & $-1/4$& $1$ &$1.4+p$&$5$&\DRAW{\PPP\NNN\PPP\NNN\PPP}&$2$&$2$&$1$ \cr
\EfiveNB & $3$ & $-1/4$& $1$ &$1.4+p$&$7$&\DRAW{\NNN\nnn\NNN\PPP\NNN\nnn\NNN}&$3/2$&$2$&$1$ \cr
\EfiveNC & $3$ & $-1/4$& $1$ &$1.4+p$&$9$&\DRAW{\PPP\ppp\PPP\ppp\PPP\ppp\PPP\ppp\PPP}&$2$&$2$&$1$ \cr
\EfiveND & $3$ & $-1/4$& $1$ &$1.4+p$&$7$&\DRAW{\NNN\PPP\ppp\PPP\ppp\PPP\NNN}&$2$ &$2$&$1$\cr
\OtwoNA & $3$ & $-1/2$& $0$ &$1.5+p$&$2$&\DRAW{\NNN\PPP}&$2$&$3$&$2$ \cr
\OfourNA & $3$ & $-1/2$& $0$ &$1.5+p$&$4$&\DRAW{\NNN\PPP\NNN\PPP}&$3/2$&$3$&$2$ \cr
\OfourNB & $3$ & $-1/2$& $0$ &$1.5+p$&$6$&\DRAW{\PPP\ppp\PPP\NNN\nnn\NNN}&$3/2$&$3$&$2$\cr
\OsixNA & $3$ & $-1/2$& $0$ &$1.5+p$&$6$&\DRAW{\PPP\NNN\PPP\NNN\PPP\NNN}&$9/8$&$3$&$2$ \crli}
\medskip
\centerline{\eightrm Table 1. Parameter values and properties of solutions.}
\bigskip
%%%%%%%%%%%%%%%%%%%%%%%%%%%%%%%%%%%%%%%%%%%%%%%%%%%
%%%%%%%%%%%%%%%%%%%%%%%%%%%%%%%%%%%%%%%%%%%%%%%%%%%

Our approach to solving the equation \equ(travelTwo)
is to turn it into a suitable fixed point problem.
This is a common strategy in computer-assisted proofs.
Alternatively, one could try a dynamical systems approach.
Traveling waves can often be viewed
as homoclinic (or heteroclinic) orbits of a dynamical system.
This approach has been successful with systems
of ordinary differential equations [\rGGKKMO,\rAKii,\rCZg,\rMatsue].
Integration methods have been applied also
to dissipative systems in infinite dimensions [\rZg,\rAKi,\rWZg]
and to some delay equations [\rMiNa,\rSZg,\rKLM,\rKiLe].
But for an advance-delay equation like \equ(travelOne),
this does not seem to be a workable approach.
In this context, we should mention that the system
\equ(FPU) is Hamiltonian.

A fixed point equation $G(v)=v$ for the function $v=\DD u$ can be obtained by
integrating both sides of \equ(travelTwo) twice.
The transformation $G$ improves regularity.
But due to the non-compactness of the domain $\real$,
we found it difficult to come up with an expansion for $v$
that allows for accurate approximations and is suitable
for a rigorous computer-assisted analysis.
The following turns out to work extremely well.
A function $v$ in one of our spaces $\Brr$
is given by a sequence of ``arcs'' $v_j=v(\bdot-j)$, indexed by integers $j$
and defined on the interval $[-\shalf,\shalf]$.
Each arc $v_j$ is real analytic on this interval
and represented by a rapidly converging Legendre series.
So $v$ is real analytic outside the set $\integer+\shalf$.
Real analyticity on all of $\real$ is obtained
if $v$ is a fixed point of $G$,
due to the regularity-improving property of the transformation $G$.

%%%%%%%%%%%%%%%%%%%%%%%%%%%%%%%%%%%%%%%%%%%%%%%%%%%%%%%%%%%%%
\vskip 0.6cm
\hbox{\hskip0.0cm
\includegraphics[height=3.8cm,width=6.5cm]{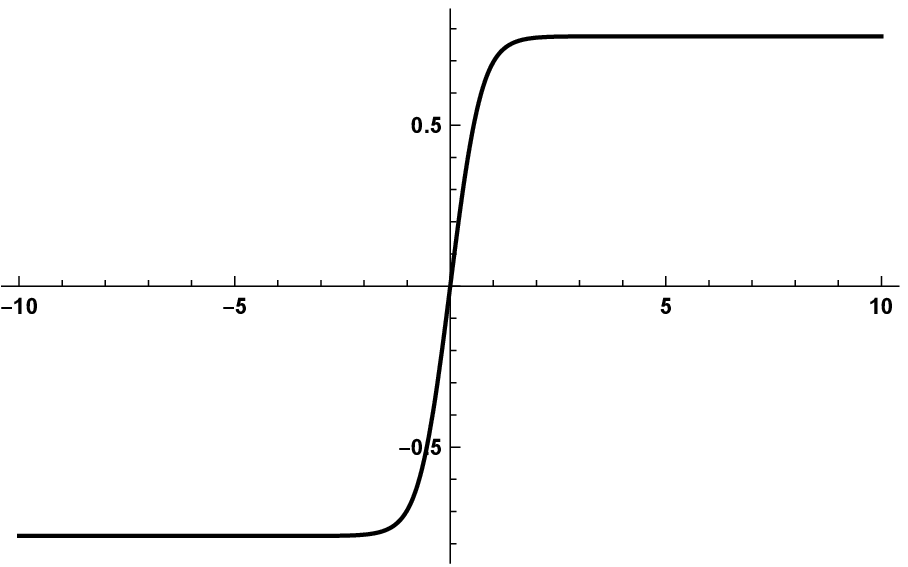}
\includegraphics[height=3.8cm,width=6.5cm]{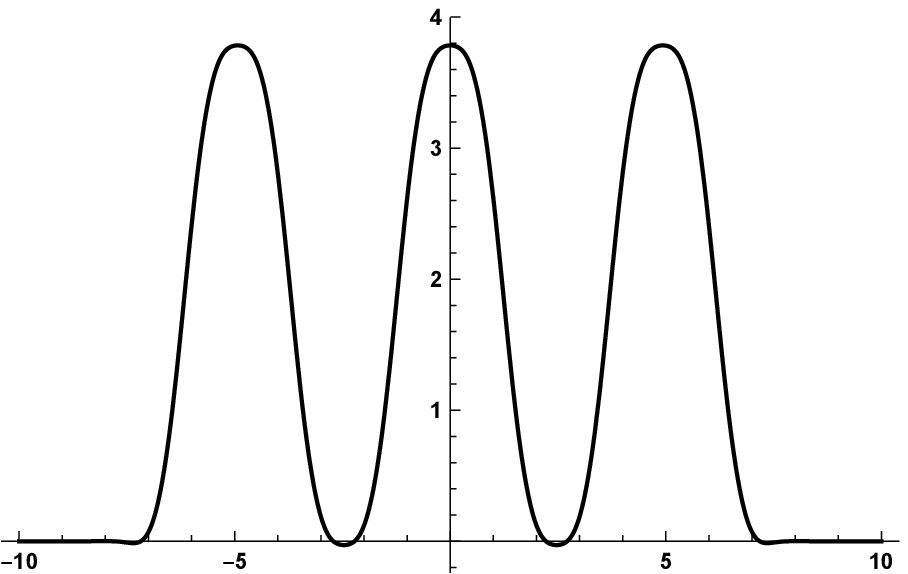}
}
\vskip 0.2cm
\leftline{\hskip 1.5cm\eightpoint{\bf Figure 1.}
Profile $u$ for the solutions \EonePA ~(left) and \OsixNA ~(right).}
%%%%%%%%%%%%%%%%%%%%%%%%%%%%%%%%%%%%%%%%%%%%%%%%%%%%%%%%%%%%%

\subsection Numerical observations
%%%%%%%%%%%%%%%%%%%%%%%%%%%%%%%%%%

As mentioned above, we solve the equation \equ(travelTwo)
by converting it to a suitable fixed point equation $G(v)=v$.
Then, starting with an approximate numerical solution $\bar v$,
we prove that $G$ has a true fixed point $v$ nearby.

To find approximate solutions we use a Newton-type iteration.
Some approximate solutions were found by
starting the iteration with a randomly generated function $v$.
For others, we started with an initial guess, produced by
combining one-bump or two-bump solutions into functions with multiple bumps.
(In this discussion, a ``bump'' is a local maximum above $\sfrac{1}{8}$
or local minimum below $-\sfrac{1}{8}$.)
Figures 4 and 6 suggest that solutions
with an arbitrary large number of bumps can be obtained this way.
To be more precise, multi-bump solutions were found only
for negative value of $\mu$.
Furthermore, the precision of our method deteriorates as $|\mu|$ approaches $1$.
But for each of the numerical solutions that we found,
it was possible to ``continue'' the solution
to other values of $\mu$ of the same sign,
without encountering any bifurcations.
So the choice of $\mu$-values in Table 1 is somewhat random,
except for the sign.

Let us call two functions $v$ and $w$ ``independent''
if they have disjoint supports, separated by an interval of length $1$.
In this case $\DD^2\phi'(v+w)=\DD^2\phi'(v)+\DD^2\phi'(w)$.
So if an approximate solution $\bar v$
is small (in modulus) on an interval of length $1$
but has bumps on both sides of this interval,
then it is a sum of two nearly independent approximate solutions.
In such cases, our Newton iteration
is pushing the two parts farther and farther apart.
This suggests that there exists no true solution with nearly independent parts.
(Near-independence is accompanied
by eigenvalues close to $1$ for the linearized problem,
so it was necessary to increase the numerical precision in such cases.)

For negative values of $\mu$, it appears that multi-bump solutions
exist only for certain specific arrangements of the bumps.
And for nonnegative values of $\mu$, no multi-bump solutions appear to exist.
As $\mu$ approaches zero from below,
the bumps approach near-independence in the sense described above;
see also Figure 8 (left).
In the case of solutions that have negative local maxima
or positive local minima,
these extreme values approach zero as $\mu\to 0$.
This behavior is illustrated in Figure 8 (right).

Solutions for $\mu<0$ appear not to change significantly
as $\mu$ approaches the value $-1$.
But the one-bump solution for $\mu>0$ widens,
and its amplitude decreases, as $\mu$ approaches $1$.
This behavior is illustrated in Figure 7.
Presumably, the function $q$ associated
with our one-bump solution $v$ approaches a small
traveling wave of the type considered in [\rIosK,\rIos] for $\mu\simeq 1$.

%%%%%%%%%%%%%%%%%%%%%%%%%%%%%%%%%
%%%%%%%%%%%%%%%%%%%%%%%%%%%%%%%%%
\section Arcs and Legendre series
%%%%%%%%%%%%%%%%%%%%%%%%%%%%%%%%%
%%%%%%%%%%%%%%%%%%%%%%%%%%%%%%%%%

In this section we describe our decomposition
of a function $w\in\rmL^2(\real)$ into local ``modes''.
In order to motivate our choices,
let us write \equ(travelTwo) as the fixed point equation
$$
v=\tau^2 A^2\phi'(v)\,,\qquad A=\DD D^{-1}\,.
\equation(travelThree)
$$
Here $D^{-1}$ is defined via integration
from some arbitrary point $x_0\in\real$.
The choice of $x_0$ does not matter,
since for any antiderivative $V$ of $v$,
$$
Av=\DD V=V\bigl(\bdot+\thalf\bigr)-V\bigl(\bdot-\thalf\bigr)
=\int_{-1/2}^{1/2} v(\bdot+s)\,ds\,.
\equation(AvConv)
$$
Notice that $A$ is a convolution operator:
If $\chi$ denotes the indicator function
of the interval $[-\shalf,\shalf]$, then $Av=\chi\ast v$.
Thus, the operator $A^2$ that appears in \equ(travelThree)
is convolution with the function $\chi\ast\chi$,
also known as the cardinal b-spline of order $2$,
with separation $1$ between the knots.
This suggests that we decompose a function $w$ on $\real$
as follows:
$$
w=\sum_{j\in\oldinteger}\proj_jw\,,\qquad
(\proj_jw)(x)=\cases{w(x), &if $x\in I_j$;\cr 0, &otherwise;\cr}
\equation(PjvDef)
$$
where $I_j=[\tinyskip j-\shalf,j+\shalf]$ for every integer $j$.

Our goal is to compute $A^2v$ via antiderivatives,
but to keep the computation as local as possible.
To see how this can be achieved,
consider a function $v$ on $\real$ that is supported on $I_0$.
If $v$ is orthogonal to the constant function on $I_0$,
then $v$ has an antiderivative $D^{-1}v$ on $\real$ that is again
supported in $I_0$.
More generally, if $v$ is orthogonal to all polynomials
of degree less than $n$,
then $v$ possesses antiderivatives $D^{-m}v$ of order $m\le n$
that are all supported on $I_0$.
This motivates the following choices.

Denote by $T$ translation by $1$; that is, $(Tw)(x)=w(x-1)$.
We will represent a function $w\in\rmL^2(\real)$ by the sequence
of ``arcs'' $w_j=T^{-j}\proj_j w$ indexed by integers $j$.
Each arc $w_j$ is supported in $I_0$,
and when regarded as a function in $\HH=\rmL^2(I_0)$,
it admits a unique expansion
$$
w_j=\sum_{n\in\oldnatural}w_{j,n}\PP_n\,,\qquad
\PP_n(x)=P_n(2x)\,,\qquad x\in I_0\,.
\equation(PnExpansion)
$$
Here $P_n$ denotes the Legendre polynomial of degree $n$.
The sequence of polynomials $n\mapsto\PP_n$
can be obtained from the sequence of monomials $n\mapsto(x\mapsto x^n)$
via the Gram-Schmidt orthogonalization process,
and then normalizing $\PP_n(\shalf)=1$.
So the scaled Legendre polynomials $\PP_0,\PP_1,\PP_2,\ldots$ constitute
a complete orthogonal set in $\HH$.
More specifically, we have
$$
\langle\PP_m,\PP_n\rangle_{\sss\HH}
={2\over 2n+1}\delta_{m,n}\,,\qquad
\langle g,f\rangle_{\sss\HH}=2\int_{-1/2}^{1/2}\ov{g(x)}f(x)\,dx\,.
\equation(SLOrtho)
$$

Two other well-known facts are the following.
The values of $\PP_n$ on the interval $I_0$ are bounded in modulus by $1$.
Furthermore, $(4n+2)\PP_n=\PP_{n+1}'-\PP_{n-1}'$ for all positive integers $n$.
Integrating both sides of this identity yields
$$
D^{-1}\PP_n={1\over 4n+2}\bigl(\PP_{n+1}-\PP_{n-1}\bigr)\,,
\qquad n\ge 1\,,
\equation(SLAntiDer)
$$
up to an additive constant.
Notice that $\PP_{n+1}$ and $\PP_{n-1}$ agree at $\pm\shalf$.
Thus, the right hand side of \equ(SLAntiDer) vanishes at $\pm\shalf$.
So we define $D^{-1}\PP_n$ for $n\ge 1$ via integration from $\shalf$ or $-\shalf$.

Applying the identity \equ(SLAntiDer) twice yields
$$
D^{-2}\PP_n=C_n^{+}\PP_{n+2}+C_n^{-}\PP_{n-2}-(C_n^{+}+C_n^{-})\PP_n\,,\qquad
n\ge 2\,,
\equation(SLDerDer)
$$
where
$$
C_n^{+}={1/4\over(2n+1)(2n+3)}\,,\qquad
C_n^{-}={1/4\over(2n+1)(2n-1)}\,.
\equation(Cnpmo)
$$
Notice that the function $D^{-2}\PP_n$ and its first derivative
vanish on the boundary of the interval $I_0$, for all $n\ge 2$.
As mentioned above,
this is one of our main reasons for having chosen
scaled Legendre polynomials for our expansion of arcs.

%%%%%%%%%%%%%%%%%%%%%%%%%%%%%%%%%%%%%%%%%%%%%%%%%%%%
%%%%%%%%%%%%%%%%%%%%%%%%%%%%%%%%%%%%%%%%%%%%%%%%%%%%
\section The operator $\boldmath A^2$ in more detail
%%%%%%%%%%%%%%%%%%%%%%%%%%%%%%%%%%%%%%%%%%%%%%%%%%%%
%%%%%%%%%%%%%%%%%%%%%%%%%%%%%%%%%%%%%%%%%%%%%%%%%%%%

The goal here is to give an explicit description
of how the operator $A^2$ acts on each term in the
decomposition \equ(PnExpansion).
To simplify notation, we define $\PP_n(x)=0$ for $|x|>\shalf$.

But first, let us consider a slight generalization
of the approach described in the previous subsection.
For a traveling wave $v$ that varies rapidly,
it can be advantageous to partition $\real$
into subintervals of length $1/\Scale$, for some integer $\Scale>1$.
Equivalently, we can reformulate our fixed point problem
in terms of the scaled function $w=v(\bdot/\Scale)$.
The resulting equation for $w$ is $w=\tau^2 A_\Scale^2\phi'(w)$, where
$$
A_\Scale=\Scale^{-1}\bigl[T^{-\Scale/2}-T^{\Scale/2}\bigr]D^{-1}\,,\qquad
A_\Scale^2=\Scale^{-2}\bigl[\tinyskip T^{-\Scale}+T^\Scale-2\id\tinyskip\bigr]D^{-2}\,.
\equation(LPAAf)
$$
Here $\Scale$ can be any positive integer. Notice that $A_1=A$.

By translation invariance and linearity,
it suffices to compute $A_\Scale^2w$ for a function $w$ of the form
$$
w(x)=\cases{(2x)^n, &if $-\half<x<\half$;\cr 0, &otherwise.}
\equation(wx)
$$
We are interested mainly in the cases $n=0$ and $n=1$,
where the identity \equ(SLDerDer) does not apply.
Let $D^{-1}w$ be the antiderivative of $w$
that vanishes at the origin,
and let $D^{-2}w$ be the antiderivative of $D^{-1}w$
that vanishes at the origin.
Then
$$
\bigl(D^{-2}w\bigr)(x)
=\cases{{1\over4(n+1)(n+2)}(2x)^{n+2}, &if $-\half<x<\half$;\cr
{1\over4(n+1)(n+2)}+{1\over2(n+1)}\bigl(x-\half\bigr), &if $x\ge\half$;\cr
(-1)^n{1\over4(n+1)(n+2)}+(-1)^n{-1\over2(n+1)}\bigl(x+\half\bigr), &if $x\le-\half$.\cr}
\equation(invDDwx)
$$
Notice that the function $D^{-2}w$ is affine
to the left and to the right of $I_0$.
When we apply $T^{-\Scale}+T^\Scale-2\id$, these affine parts of $D^{-2}w$
cancels at a distance $\Scale$ or larger from $I_0$.
Applying $T^{-\Scale}+T^\Scale-2\id$ to the part of $D^{-2}w$
that is supported in $I_0$ results in three copies:
one in $I_0$, with a factor $-2$,
and a translated copy in each of the intervals $I_{\pm\Scale}$.

In what follows, we restrict to $\Scale=1$ or $\Scale=2$.
After a trivial but tedious computation,
we end up with the following expressions.
For $n<2$ and $\Scale=1$ we have
$$
\eqalign{
A_1^2\PP_0
&=\textstyle
T^{-1}\bigl[{1\over 12}\PP_2+{1\over 4}\PP_1+{1\over 6}\PP_0\bigr]
+\bigl[-{1\over 6}\PP_2+{2\over 3}\PP_0\bigr]\cr
&\quad\textstyle
+T\bigl[{1\over 12}\PP_2-{1\over 4}\PP_1+{1\over 6}\PP_0\bigr]\,,\cr
A_1^2\PP_1
&=\textstyle
T^{-1}\bigl[{1\over 60}\PP_3-{1\over 10}\PP_1-{1\over 12}\PP_0\bigr]
+\bigl[-{1\over 30}\PP_3+{1\over 5}\PP_1\bigr]\cr
&\quad\textstyle
+T\bigl[{1\over 60}\PP_3-{1\over 10}\PP_1+{1\over 12}\PP_0\bigr]\,.\cr}
\equation(AAPoiOne)
$$
For $n<2$ and $\Scale=2$ we obtain
$$
\eqalign{
4A_2^2\PP_0
&=T^{-2}\Bigl[\textstyle{1\over 12}\PP_2+{1\over 4}\PP_1+{1\over 6}\PP_0\Bigr]
+T^{-1}\Bigl[\textstyle{1\over 2}\PP_1+\PP_0\Bigr]\cr
&\quad+\Bigl[\textstyle-{1\over 6}\PP_2+{5\over 3}\PP_0 \Bigr]
+T\Bigl[\textstyle-{1\over 2}\PP_1+\PP_0 \Bigr]
+T^2\Bigl[\textstyle {1\over 12}\PP_2-{1\over 4}\PP_1+{1\over 6}\PP_0\Bigr]\,,\cr
4A_2^2\PP_1
&=T^{-2}\Bigl[\textstyle{1\over 60}\PP_3-{1\over 10}\PP_1-{1\over 12}\PP_0\Bigr]
+T^{-1}\Bigl[\textstyle-{1\over 6}\PP_0\Bigr]\cr
&\quad+\Bigl[\textstyle-{1\over 30}\PP_3+{1\over 5}\PP_1\Bigr]
+T\Bigl[\textstyle{1\over 6}\PP_0\Bigr]
+T^2\Bigl[\textstyle{1\over 60}\PP_3-{1\over 10}\PP_1+{1\over 12}\PP_0\Bigr]\,.\cr}
\equation(AAPoiTwo)
$$
For $n\ge 2$ it suffices to combine the identities \equ(LPAAf) and \equ(SLDerDer).

\medskip
In our computations [\rFiles] we apply these formulas to finitely many terms in the expansion
$w=\sum_{j,n} w_{j,n}T^j\PP_n$.
In order to estimate truncation errors, we use (among others)
the following simple facts.
Consider a convolution operator of the form
$$
(Bw)_{j,n}=\sum_{i,m}B^{n,m}_{j-i}w_{i,m}\,,
\equation(Bwjn)
$$
where $i,j,m,n$ denote integers, with $m,n$ positive.
Since we are interested in upper bounds only,
assume for simplicity that the numbers $B^{n,m}_k$ and $w_{i,m}$ are all nonnegative.
We consider bounds of the type
$$
W_{J,n}(w)=\sum_{j\ge J}r^{j-J}w_{j,n}
\equation(WJn)
$$
with $r>1$, and weighted sums (over $n$) of such bounds.

\claim Proposition(WBw) $W_{J,n}(Bw)=(BW(w))_{J,n}$ for all $J$ and $n$.
If $w_{j,n}=0$ for all $j<L$, then $W_{J,n}(w)\le r^{L-J}W_{L,n}$ for all $J$;
and equality holds for $J\le L$.

A proof of these properties is straightforward.
Notice that $w\mapsto W(w)$ is a convolution,
which explains most of why $WB=BW$.

%%%%%%%%%%%%%%%%%%%%
%%%%%%%%%%%%%%%%%%%%
\section Analyticity
%%%%%%%%%%%%%%%%%%%%
%%%%%%%%%%%%%%%%%%%%

Here we describe the function spaces that are used
in our analysis of the fixed point equation \equ(travelThree).
Given $\rho\ge 1$, denote by $\AA_\rho$ the space of all functions
$f$ in $\HH=\rmL^2(I_0)$ that have a finite norm
$$
\|f\|_\rho\defeq\sum_{n\in\oldnatural}|f_n|\rho^n\,,\qquad
f=\sum_{n\in\oldnatural}f_n\PP_n\,.
\equation(AANorm)
$$
The even and odd subspaces of $\AA_\rho$
are denoted by $\AA_\rho^\even$ and $\AA_\rho^\odd$, respectively.
If $\rho>1$, then every function in $\AA_\rho$
admits an analytic continuation to
the complex open neighborhood $\EE_\rho$ of $I_0$
whose boundary is the ellipse $\bigl|4z+z^{-1}\bigr|=4\rho$.
This is a consequence of Bernstein's lemma.

Consider now functions $w:\real\to\real$ that admit an expansion \equ(PjvDef)
with $w_j\in\AA_\rho$ for all $j$.
Given $r\ge 1$, define $\Brr$ to be the space
of all such functions $w$ that have a finite norm
$$
\|w\|_{\rho,r}\defeq\sum_{j\in\oldinteger}\|w_j\|_\rho r^{|j|}
=\sum_{j\in\oldinteger}\;\sum_{n\in\oldnatural}|w_{j,n}|r^{|j|}\rho^n\,.
\equation(BBNorm)
$$
The even and odd subspaces of $\Brr$
are denoted by $\Brr^{\,\even}$ and $\Brr^{\,\odd}$, respectively.
Since our fixed point equation \equ(travelThree)
involves a product of functions, we will need the following.

\claim Proposition(BanachAlg)
The spaces $\AA_\rho$, $\AA_\rho^\even$, $\Brr$, and $\Brr^{\,\even}$
are Banach algebras.

\proof
The product of two function in $\AA_\rho$
can be estimated by using the linearization formula
$$
\PP_k\PP_l=\sum_m C_{k,l,m}\PP_m\,,\qquad
\sum_m C_{k,l,m}=1\,.
\equation(LegendreLinearization)
$$
The existence of such an expansion follows by completeness,
and the second identity is a consequence of the fact that $\PP_n(\shalf)=1$ for all $n$.
Clearly $C_{k,l,m}=0$ for $m>k+l$.
And by parity, we have $C_{k,l,m}=0$ unless $k+l+m$ is even.
Below will see that the coefficients $C_{k,l,m}$ are all nonnegative.
So if $f$ and $g$ belong to $\AA_\rho$, then
$$
\eqalign{
\|fg\|_\rho
&=\sum_{m\in\oldnatural}\,\Biggl|\,\sum_{k,l\in\oldnatural}
C_{k,l,m}\bigl(\rho^k f_k\bigr)\bigl(\rho^l g_l\bigr)\rho^{m-k-l}\Biggr|\cr
&\le\sum_{k,l\in\oldnatural}\rho^k|f_k|\rho^l|g_l|\sum_{m\le k+l}|C_{k,l,m}|
=\|f\|_\rho\|g\|_\rho\,.\cr}
\equation(AABanachAlg)
$$
An analogous inequality $\|vw\|_{\rho,r}\le\|v\|_{\rho,r}\|w\|_{\rho,r}$
for functions $v,w\in\Brr$ follows trivially.

\smallskip
The following expression for the coefficients $C_{k,l,m}$
was first found by Adams [\rAdams].
$$
C_{k,l,m}={a(s-k)a(s-l)a(s-m)\over a(s)}\,{2m+1\over 2s+1}\,,
\qquad a(n)=2^{-n}{2n\choose n}\,,
\equation(CklmAdams)
$$
where $s=\half(k+l+m)$.
Our main reason for giving this formula here
is that it is being used in our programs [\rFiles].
But it also shows that the coefficients $C_{k,l,m}$ are all nonnegative,
as claimed above.
\qed

\demo Remark(Eval)
A noteworthy consequence of \equ(AABanachAlg) is the following.
Let $f\in\AA_\rho$. Then the spectral radius of the operator
$g\mapsto fg$ acting on $\AA_\rho$ is bounded by $\|f\|_\rho$.
This implies e.g.~that the analytic continuation of $f$
satisfies $|f(z)|\le\|f\|_\rho$ for all $z\in\EE_\rho$.

Notice that the translation operators $T$ and $T^{-1}$ are bounded on $\Brr$.
So the following is essentially a consequence
of the identity \equ(SLDerDer) and \equ(Cnpmo).

\claim Proposition(AAbounded)
The operator $A_\Scale$ defined by \equ(LPAAf) is bounded on $\Brr$.

The following will be used to prove
that a solution $v\in\Brr$ of the equation \equ(travelTwo)
is real analytic.

\claim Proposition(BBanalytic)
Assume that $\rho>1$.
Then there exists $\epsilon>0$ such that,
if $w\in\Brr$ is of class $\rmC^\infty$,
then $w$ extends analytically to the strip $|\Im z|<\epsilon$.

\proof
Let $w\in\Brr$.
As described after \equ(AANorm), each arc $w_j$ extends
analytically to a domain $\EE_\rho$ that includes an open rectangle
$\RR=(-\shalf-\epsilon,\shalf+\epsilon)\times(-\epsilon,\epsilon)$.
Thus $\proj_j v$ extends analytically to $T^j\RR$, for each $j\in\integer$.
Denote this extension by $E_j w$.
Then the domains of both $E_j w$ and $E_{j+1}w$ include the rectangle
$\RR_j=(j+\shalf-\epsilon,j+\shalf+\epsilon)\times(-\epsilon,\epsilon)$.
Assume now that $w$ belongs to $\rmC^\infty(\real)$.
Then the derivatives $D^nE_j w$ and $D^nE_{j+1}w$
agree at the point $j+\shalf$, for all $n\ge 0$.
Thus, $E_j w$ and $E_{j+1}w$ agree on $\RR_j$.
This holds for all $j\in\integer$,
so $w$ extends analytically to $\real\times(-\epsilon,\epsilon)$.
\qed

%%%%%%%%%%%%%%%%%%%%%%%%%%%%%%%%%
%%%%%%%%%%%%%%%%%%%%%%%%%%%%%%%%%
\section The fixed point equation
%%%%%%%%%%%%%%%%%%%%%%%%%%%%%%%%%
%%%%%%%%%%%%%%%%%%%%%%%%%%%%%%%%%

Consider the parameter values $(m,\bar\mu,\sigma,r,k,\Scale)$
associated with some fixed row in Table 1.
In addition, we choose $\rho=\sfrac{17}{16}$.
Let $\mu$ be a real number of modulus $|\mu|<1$.
Then the equation \equ(travelThree) for these parameters
reads $v=\nu A^2v^m+\mu A^2v$.
The corresponding equation for the function $w=v(\bdot/\Scale)$ is
$w=\nu A_\Scale^2w^m+\mu A_\Scale^2w$, with $A_\Scale$ as defined by \equ(LPAAf).
This equation for $w$ is equivalent to the fixed point equation $G_q(w)=w$,
where
$$
G_q(w)=\nu A_\Scale^2\Sigma_q w^m+\mu^k A_\Scale^{2k}w\,,\qquad
\Sigma_q=\sum_{n=0}^{k-1}\mu^n A_\Scale^{2n}\,.
\equation(GqDef)
$$
The subscript used here is $q=(m,\mu,k,\Scale)$.
By Propositions \clmno(BanachAlg) and \clmno(AAbounded),
$G_q$ is differentiable as a map on $\Brr$.
Notice that, if $w$ and $h$ belong to $\Brr$,
and if $h$ is supported far away from the origin,
then $DG_q(w)h$ is approximately equal to $\mu^k A_\Scale^{2k}h$.
So the transformation $G_q$ contracts tails by roughly a factor $|\mu|^k$.
This should make clear why we need larger values of $k$
when $|\mu|$ is close to $1$.

\claim Definition(SplinesEtc)
A compactly supported function $w:\real\to\real$
whose arcs $w_j$ are all polynomials will be called a spline.

Let $\tau=1-\sigma$.
Our choice of the exponent $m$ in \equ(muvnuvm) guarantees
that $G_q(w)$ has parity $\tau$ whenever $w$ has parity $\tau$.
Given a function $\bar w\in\Brrt$,
and a bounded linear operator $M$ on $\Brrt$, define
$$
\NN_q(h)=G_q(\bar w+\Lambda h)-\bar w+Mh\,,\qquad\Lambda=\id-M\,,
\equation(contraction)
$$
for every function $h\in\Brrt$.
Clearly, if $h$ is a fixed point of $\NN_q$
then $\bar w+\Lambda h$ is a fixed point of $G_q$.
For practical reasons, we choose $\bar w$ to be a spline
that is an approximate fixed point of $G_{\bar q}$
with $\bar q=(m,\bar\mu,k,\Scale)$.
And for $M$ we choose a finite rank operator
such that $\Lambda$ is an approximate inverse of $\id-DG_{\bar q}(\bar w)$.
If $\mu$ is sufficiently close (but not necessarily equal) to $\bar\mu$,
then we can expect $\NN_q$ to be a contraction near the origin.

\claim Lemma(contr)
Consider a fixed but arbitrary row in Table 1
and the parameter values given in that row.
Let $\rho=\sfrac{17}{16}$ and $\tau=1-\sigma$.
Then there exists a spline $\bar w$,
a bounded linear operator $M$ on $\Brrt$,
and positive constants $\eps,K,\delta$ satisfying
$\eps+K\delta<\delta$, such that for every value of $\mu$
in some open neighborhood of $\bar\mu$,
the transformation $\NN_q$ defined by \equ(contraction) satisfies
$$
\|\NN_q(0)\|_{\rho,r}\le\eps\,,\qquad
\|D\NN_q(h)\|_{\rho,r}\le K\,,\qquad h\in B_\delta\,,
\equation(contr)
$$
where $B_\delta$ denotes the closed ball of radius $\delta$
in $\Brrt$, centered at the origin.
Furthermore, if $w=\bar w+\Lambda h$ with $h\in B_\delta$,
then the function $v=w(\Scale\tinyskip\bdot)$ has the properties
listed in (the given row of) Table 1 concerning
the $\sup$-norm and the local extrema.

Notice that $\|w-\bar w\|_{\rho,r}\le\delta'$, where $\delta'=\|\Lambda\|\delta$.
Our proof of \clm(contr) yields $\delta'<2^{-32}$ for all solutions.
This bound can be made as small as desired
by running our programs (which also determine $\bar w$)
at higher numerical precision.

Based on this lemma, we can now give a

\proofof(main)
By the contraction mapping principle,
the given bounds imply that $\NN_q$
has a unique fixed point $h$ in $B_\delta$.
The corresponding function $w=\bar w+\Lambda h$ is a fixed point of $G_q$.
Given that $A_\Scale^2$ includes two antiderivatives,
the identity $w=G_q(w)$ implies that $w$ is of class $\rmC^\infty$.
So by \clm(BBanalytic), $w$ extends to an analytic function on some strip $|\Im z|<\epsilon$.
This extension still decreases exponentially:
by \dem(Eval) we have a uniform bound $|w(x+iy)|\le Cr^{-|x|}$
that holds for all $x,y\in\real$ with $|y|<\epsilon$.
These properties of $w$ imply that the function $v=w(\Scale\tinyskip\bdot)$
extends analytically to the strip $|\Im z|<\epsilon/\Scale$,
decreases exponentially, and satisfies the the equation \equ(travelTwo).

Consider now the function $u=u_\ssL$ defined by the equation \equ(uLuR).
The above-mentioned properties of $v$ imply that $u$ is real analytic.
Furthermore, $\DD u=v$.
The equation \equ(travelTwo) implies that
the function $g$ defined by $g=u''-\tau^2\DD\phi'(v)$ satisfies $\DD g=0$.
Thus $g$ is periodic with period $1$.
The function $\DD\phi'(v)$ vanishes at $\pm\infty$,
so $u''$ approaches the periodic function $g$ at $\pm\infty$.
But $u''=u_\ssL''$ vanishes at $-\infty$.
Thus $g=0$, which implies that $u''=\tau^2\DD\phi'(\DD u)$.
In other words, $u$ satisfies the equation \equ(FPU).

The same arguments apply to the function $u=u_\ssR$.
The difference $f=u_\ssL-u_\ssR$ is periodic, since $\DD f=0$.
But $f''$ vanishes at infinity, so $f$ is constant.
Now define $u=\shalf u_\ssL+\shalf u_\ssR$.
Then $u$ is real analytic, satisfies the equation \equ(FPU),
and has parity $\sigma$.

The claims in \clm(main)
concerning the $\sup$-norm and the local extrema
of the function $v=\DD u$ follow from the last statement in \clm(contr).
\qed

%%%%%%%%%%%%%%%%%%%%%%%%%%%
%%%%%%%%%%%%%%%%%%%%%%%%%%%
\section Computer estimates
%%%%%%%%%%%%%%%%%%%%%%%%%%%
%%%%%%%%%%%%%%%%%%%%%%%%%%%

What remains to be done is to verify \clm(contr).
This is carried out with the aid of a computer.
To be more specific, consider the parameter values $(m,\bar\mu,\sigma,r,k,\Scale)$
from a fixed but arbitrary row in Table 1.
As a first step, we determine an approximate fixed point $\bar w$ of $G_q$
and an approximate inverse of $\id-DG_{\bar q}(\bar w)$
of the form $\Lambda=\id-M$, with $M$ of finite rank.
The remaining steps are rigorous:
We compute an upper bound $\eps$ on the norm of $\NN_q(0)$,
and an upper bound $K$ on the operator norm of $D\NN_q(h)$
that holds for all $h$ of norm $4\eps$ or less.
This is done simultaneously for all values of $\mu$
in some open interval centered at $\bar\mu$.
After verifying that $K<\sfrac{7}{8}$, we choose a positive $\delta<8\eps$
in such a way that $\eps+K\delta<\delta$.
The last statement in \clm(contr)
is verified by estimating $v(x_i)$ at a finite number of points $x_i\in\real$.

The rigorous part is still numerical,
but instead of truncating series and ignoring rounding errors,
it produces guaranteed enclosures at every step along the computation.
This part of the proof is written in the
programming language Ada [\rAda].
The following is meant to be a rough guide for the reader who wishes to check the
correctness of our programs.
The complete details can be found in [\rFiles].

In the present context, a ``bound'' on a map $f:\XX\to\YY$
is a function $F$ that assigns to a set $X\subset\XX$
of a given type ({\tt Xtype}) a set $Y\subset\YY$
of a given type ({\tt Ytype}), in such a way that
$y=f(x)$ belongs to $Y$ for all $x\in X$.
In Ada, such a bound $F$ can be implemented by defining
a {\tt procedure F(X\col in Xtype\scol Y\col out Ytype)}.

For balls in a real Banach algebra $\XX$ with unit ${\bf 1}$,
we use a data type {\tt Ball}.
A {\tt Ball} consists of a pair {\tt B=(B.C,B.R)},
where {\tt B.C} is a representable number ({\tt Rep})
and {\tt B.R} a nonnegative representable number ({\tt Radius}).
The corresponding ball in $\XX$ is the set
${\tt B}_{\sss\XX}=\{x\in\XX:\|x-({\tt B.C}){\bf 1}\|\le{\tt B.R}\}$.
Our bounds on some standard functions involving the type {\tt Ball}
are defined in the packages {\tt Std\_Balls}.
Other basic functions are covered in the packages {\tt Vectors} and {\tt Matrices}.
Bounds of this type have been used in many computer-assisted proofs;
so we focus here on the more problem-specific aspects of our programs.

\subsection Analytic arcs
%%%%%%%%%%%%%%%%%%%%%%%%%

Consider the space $\AA_\rho$ for a given {\tt Radius} $\rho\ge 1$.
Our enclosures for functions in $\AA_\rho$ are associated with a data type {\tt Legend},
based on scalars of type {\tt Ball}, with $\XX=\oldreal$.
Given a fixed odd integer $D>1$, a {\tt Legend} is in essence a pair {\tt G=(G.C,G.E)},
where {\tt G.C} is an {\tt array(0\dotdot D) of Ball}
and {\tt G.E} is an {\tt array(0\dotdot D+2) of Radius}.
The corresponding set ${\tt G}_{\!\sss\AA}\subset\AA_\rho$
consists of all functions $g$ that admit a representation
$$
g=\sum_{n=0}^D c_n\PP_n+\sum_{m=0}^{D+2}\rho^m E_m\,,
\equation(LegendFun)
$$
with $c_n\in{\tt G.C}(n)_{\sss\oldreal}$ and $E_m\in\AA_{\rho,m}$,
satisfying $\|E_m\|_\rho\le{\tt G.E}(m)$.
Here $\AA_{\rho,m}$ denotes the subspace of $\AA_\rho$
consisting of all functions $E\in\AA_\rho$
that are orthogonal to all polynomials of degree less than $m$
and have the same parity as $\PP_m$.

The type {\tt Legend} is defined in the package {\tt Legends}.
This package also implements basic bounds on functions to/from the space $\AA_\rho$.
This includes a bound {\tt Prod} on the product $(f,g)\mapsto fg$,
based on the identities \equ(LegendreLinearization) and \equ(CklmAdams).
These bounds are quite straightforward, so we refer to [\rFiles] for details.

In {\tt Legends.Chain} we use \equ(SLDerDer) to define a bound {\tt DDInvHigh}
on the operator $D^{-2}$ on $\AA_\rho$,
restricted to function \equ(LegExpansion) with $f_0=f_1=0$.
For compactly supported functions $w\in\AA_\rho^{\oldinteger}$
we use enclosures of a type {\tt LVector}.
This type is defined as an unconstrained {\tt array (Integer range <>) of Legend}.
Using {\tt DDInvHigh},
as well as the identities \equ(LPAAf), \equ(AAPoiOne), and \equ(AAPoiTwo),
we define a bound {\tt AA} on the linear operator $A_\Scale^2$ for such functions.

\subsection Piecewise real analytic functions
%%%%%%%%%%%%%%%%%%%%%%%%%%%%%%%%%%%%%%%%%%%%%

Let $L$ be a fixed integer larger than $1$.
Then a function $w\in\Brrt$ has a unique decomposition
$$
w=\sum_{|j|<L}T^j w_j+G\,,
\equation(LChainFun)
$$
with $G\in\Brrt$ supported outside $(\shalf-j,j-\shalf)$.
The sum in this equation will be referred to as the ``center'' of $w$,
and $G$ will be referred to as the ``tail'' of $w$.

Let $r\ge 1$ be a fixed {\tt Radius}.
An enclosure for a tail in $\Brrt$ is defined by a {\tt Legend} {\tt G}
with the property that ${\tt G.C\hbox{$(n)$}.C}$ is zero for all $n$.
The corresponding set ${\tt G}_{\sss\BB}\subset\Brrt$
consists of all functions $G$ that admit a representation
$$
G=r^L\sum_{|j|\ge L}T^j g_j\,,\qquad
g_j=\sum_{n=0}^D c_{j,n}\PP_n+\sum_{m=0}^{D+2}\rho^m E_{j,m}\,,
\equation(TailFun)
$$
with coefficients $c_{j,m}\in\real$ and functions $E_{j,m}\in\AA_{\rho,m}$
satisfying the bounds
$$
\sum_{j\ge L} r^{j-L}|c_{j,n}|\le{\tt G.C\hbox{$(n)$}.R}\,,\qquad
\sum_{j\ge L} r^{j-L}\|E_{j,m}\|_\rho\le{\tt G.E\hbox{$(m)$}}\,.
\equation(TailBound)
$$
The coefficients $c_{j,m}$ and functions $E_{j,m}$
for $j\le-L$ are determined by the requirement that $G$ has parity $\tau$.

For more general subsets of $\Brrt$ we use a data type {\tt LChain},
which consists of a triple {\tt W=(W.R,W.Par,W.C)},
where ${\tt W.R}=r$, ${\tt W.Par}=\tau$,
and where {\tt W.C} is an {\tt LVector(0\dotdot L)}.
The component {\tt W.C(0)} must have parity $\tau$.
The corresponding set ${\tt W}_{\sss\BB}$ is defined as
the set of all functions \equ(LChainFun),
where $w_j\in{\tt W.C\hbox{$(j)$}}_{\sss\AA}$ for $0\le j<L$,
and $G\in{\tt W.C\hbox{$(L)$}}_{\sss\BB}$.
The arcs $w_j$ for $-L<j<0$ are determined
by the requirement that $w$ has parity $\tau$.

These types are defined in the package {\tt Legends.Chain}
which takes ${\tt JCMax}=L-1$ and ${\tt Scale}=\Scale$ as arguments.
This package also implements basic bounds on functions to/from the spaces $\Brrt$.
Most are straightforward combinations of bounds defined in {\tt Legends},
such as the bound {\tt Norm} on $w\mapsto\|w\|_{\rho,r}$
or the bound {\tt Prod} on $(v,w)\mapsto vw$.
The representation \equ(TailFun) for the tail $G$ has been chosen in such a way
that the tail component {\tt W.C(L)} of an {\tt LChain} {\tt W} can often be treated
the same way as the other components {\tt W.C(J)}.

\subsection Transformations and their derivatives
%%%%%%%%%%%%%%%%%%%%%%%%%%%%%%%%%%%%%%%%%%%%%%%%%

Consider first the operators $A_\Scale^{2k}$ and $A_\Scale^2\Sigma_q$
that appear in the definition \equ(GqDef) of the transformation $G_q$.
Our bounds on these two operators are given
by the two procedures {\tt AAPower} and {\tt SumAAPowers}.
They are more elaborate than the bounds discussed so far,
due to the fact that $A_\Scale^2$ is nonlocal.
In particular, $A_\Scale^2$ couples the center and tail of a function $w$.

Consider the task of implementing a bound on $A_1^2$.
Given an {\tt LChain} {\tt W},
consider a fixed but arbitrary function $w\in{\tt W}_{\sss\BB}$
of the form \equ(LChainFun).
The goal is to find an enclosure ${\tt U}_{\sss\BB}$
for the function $u=A_1^2w$ that only depends on {\tt W}.
By linearity, we can consider centers and tails separately.
Assume first that {\tt W} has a zero tail {\tt W.C(L)}.
Then, by using the above-mentioned bound {\tt AA} for compactly supported chains,
we obtain a {\tt LVector}-type enclosure ${\tt P(0\dotdot L)}_{\AA^\oldinteger}$
for $u$. Setting {\tt U.C(0\dotdot L-1) := P(0\dotdot L-1)}
and converting {\tt P(L)} to a tail {\tt U.C(L)},
we obtain the desired enclosure ${\tt U}_{\sss\BB}$.
This part, generalized to $A_\Scale^{2k}$, is implemented by the procedure {\tt LA\_AAPower}.
Next, consider the case where {\tt W.C(J)} is zero for all {\tt J$<$L}.
In order to determine the center part of {\tt U},
the tail {\tt W(L)} can be considered to be a bound on $w_{\sss L}$ only,
since $A_1^2$ is a convolution with a kernel supported in $[-1,1]$.
So the center components {\tt U.C(0\dotdot L-1)} of {\tt W}
can be obtained again via the procedure {\tt AA}.
This part, generalized to $A_\Scale^{2k}$, is implemented by {\tt HL\_AAPower}.
An enclosure ${\tt U.C(L)}_{\sss\BB}$ for the tail of $u$
is constructed in the procedure {\tt HH\_AAPower}.
In this part we use \clm(WBw).
The enclosure {\tt U} for $u=A_\Scale^{2k}w$ returned by {\tt AAPower}
is the {\tt Sum} of the enclosures returned by
{\tt LA\_AAPower}, {\tt HL\_AAPower} and {\tt HH\_AAPower}.
Our bound {\tt SumAAPowers} on $A_\Scale^2\Sigma_q$ is very similar.
For more details we refer to the program code [\rFiles].

Our bounds {\tt GMap} and {\tt DGmap} on the map $G_q$ and its derivative, respectively,
are defined in the package {\tt Legends.Chain.Fix}.
They are obtained simply by combining
lower-level bounds like {\tt Prod}, {\tt AAPower}, and {\tt SumAAPowers}.
The construction \equ(contr) of a quasi-Newton map $\NN$
from a given map $G$ is sufficiently general and useful
that it has been implemented in a generic package {\tt Linear.Contr}.
The same package has been used before in [\rAKiv,\rAKv].
Our instantiation of {\tt Linear.Contr}
defines bounds {\tt Contr} and {\tt Contr}
on the transformation $\NN_q$ and its derivative, respectively.
The type {\tt LMode} that is used to instantiate {\tt Linear}
will be described below.

\subsection Operator norms
%%%%%%%%%%%%%%%%%%%%%%%%%%

Consider the task of estimating the norm of a linear operator on $\Brrt$.
Let $\hat\tau=1-2\tau$.
Denote by $\SS$ the set of all pairs $s=(j,n)$ of integers $j,n\ge 0$
with the property that $n\equiv\tau\,(\mod 2)$ whenever $j=0$.
If we define $\varrho_{(j,n)}=r^j\rho^n$ and
$$
h^{(j,n)}(t)=\cases{\thalf r^{-j}\rho^{-n}[\PP_n(t-j)+\hat\tau\tinyskip\PP_n(j-t)],
&if $t\in I_j\cup I_{-j}$;\cr 0, &otherwise;\cr}
\equation(BBbasis)
$$
then a function $w\in\Brrt$
and its norm $\|w\|=\|w\|_{\rho,r}$ can be written as
$$
w=\sum_{s\in\SS}w_sh^s\,,\qquad
\|w\|=\sum_{s\in\SS}|w_s|\varrho_s\,.
\equation(wNorm)
$$
A useful feature of such weighted $\ell^1$ spaces is the following.
Let $\LL$ be a continuous linear operator on $\Brrt$.
Then the operator norm of $\LL$ is simply $\|\LL\|=\sup_{s\in\SS}\,\bigl\|\LL h^s\bigr\|$.
In order to estimate this norm,
we first choose a suitable partition $\{S_1,S_2,\ldots,S_M\}$ of $\SS$.
Then
$$
\|\LL\|=\max\{b_1,b_2,\ldots,b_M\}\,,\qquad
b_m=\sup_{s\in S_m}\,\bigl\|\LL h^s\bigr\|\,.
\equation(LLNorm)
$$
The sets $S_m$ that we use in our partitions of $\SS$
are specified by data of type {\tt LMode}.
A partition is represented by an {\tt array (1\dotdot M) of LMode}.
Such partitions are created by the procedure {\tt Make} in {\tt Legends.Chain}.
To simplify notation, let us identify a {\tt LMode} {\tt S}
with the corresponding subset $S\subset\SS$.
A {\tt procedure Assign(S\col in LMode\scol H\col out LChain)}
defines a set ${\tt H}_{\sss\BB}\subset\Brrt$
that contains all functions $h^s$ with $s\in S$.

The way this is being used is as follows.
Let {\tt LinOp} be a bound on the operator $\LL$.
Then a {\tt Ball}-type enclosure ${\tt B}_{\sss\oldreal}$ for the constant $b_m$
is obtained by calling {\tt Assign(S,H)} with $S=S_m$,
followed by {\tt LinOp(H,G)} and then {\tt Norm(G,B)}.
For the operator $\LL=D\NN_q(w)$,
this is carried out by the procedure {\tt DContrNorm} in {\tt Legends.Chain.Fix}.
This procedure is little more than an instantiation
of the procedure {\tt Op\_Norm} from the generic package {\tt Linear},
with {\tt LinOp} being in essence {\tt DContr}.

\bigskip
For the complete details we refer to the source code of our programs [\rFiles].
For the set of representable numbers ({\tt Rep}) we
choose standard extended floating-point numbers [\rIEEE] that support controlled rounding,
and for bounds on non-elementary {\tt Rep}-operations we use the open source MPFR library
[\rMPFR]. Our programs were run successfully on a standard desktop machine, using a public
version of the gcc/gnat compiler [\rGnat].

\vfil\break

%%%%%%%%%%%%%%%%%
%%%%%%%%%%%%%%%%%
\section Appendix
%%%%%%%%%%%%%%%%%
%%%%%%%%%%%%%%%%%

The figures below show graphs of the functions $v=\DD u$
associated with our solutions $u$ of the equation \equ(travelOne).
Recall that $v$ is a solution of the equation
$$
v=A^2\bigl(\mu v+v^m\bigr)\,,\qquad A=\DD D^{-1}\,.
\equation(travelTwoAgain)
$$
Notice that odd solutions ($\sigma=0$) of this equation require $m$ to be odd.

All solutions depicted here are for the $\beta$-model ($m=3$).
The solutions \EonePZ ~and \EtwoNZ ~for the $\alpha$-model ($m=2$)
are similar to the $\beta$-model solutions \EonePA ~and \EtwoNA.

\subsection Some solutions with $\boldmath\sigma=1$
%%%%%%%%%%%%%%%%%%%%%%%%%%%%%%%%%%%%%%%%%%%%%%%%%%%

%%%%%%%%%%%%%%%%%%%%%%%%%%%%%%%%%%%%%%%%%%%%%%%%%%%%%%%%%%%%%
\vskip0.2cm
\hbox{\hskip0.0cm
\includegraphics[height=3.2cm,width=2.4cm]{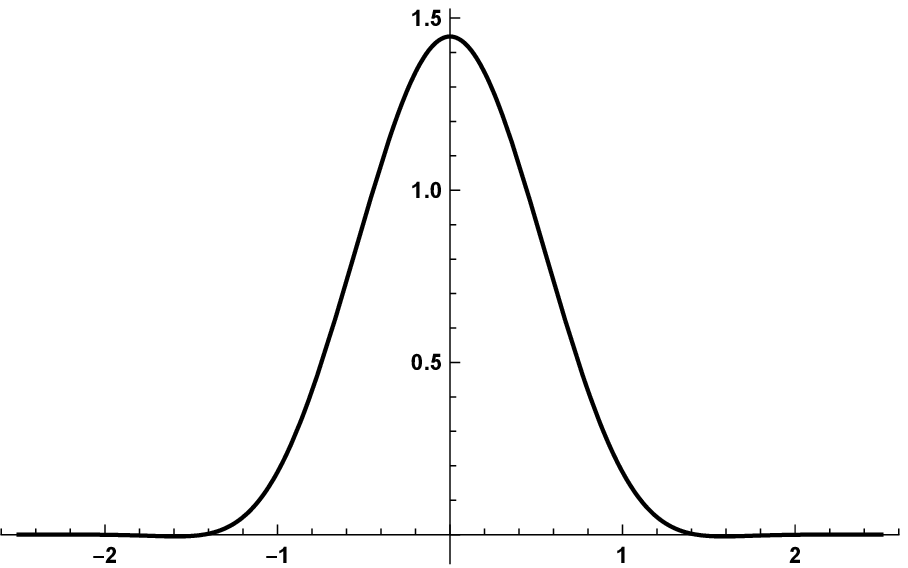} %% mu=-1/4
\includegraphics[height=3.2cm,width=2.4cm]{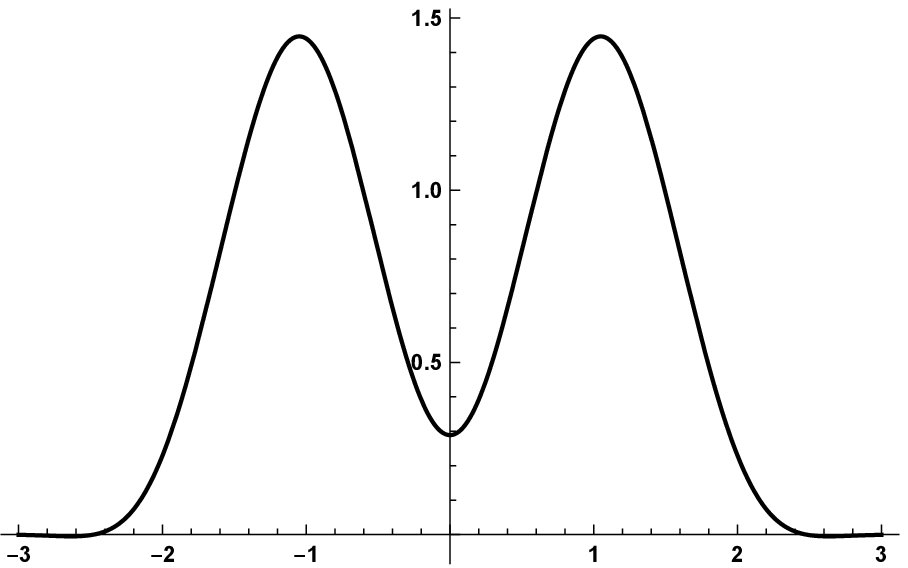} %% mu=-1/4
\includegraphics[height=3.2cm,width=2.4cm]{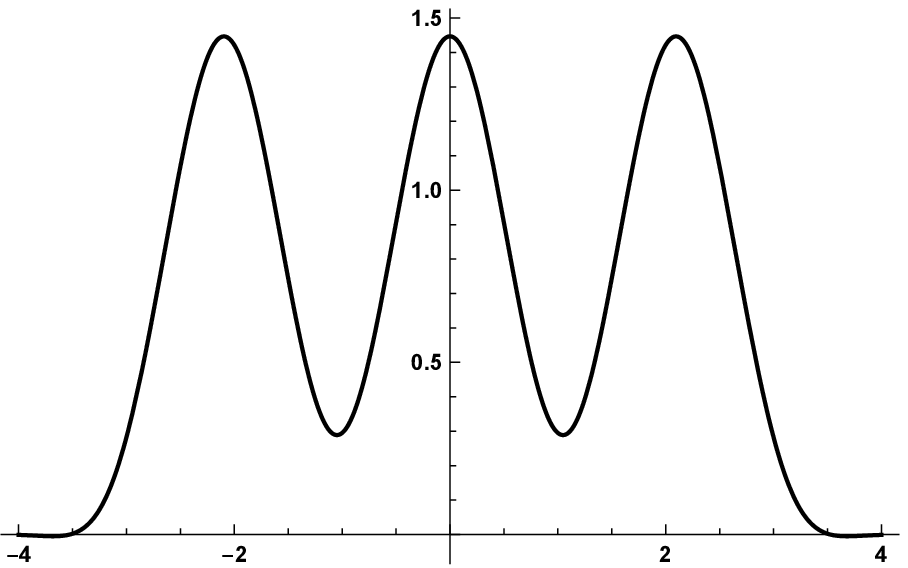} %% mu=-1/4
\includegraphics[height=3.2cm,width=2.6cm]{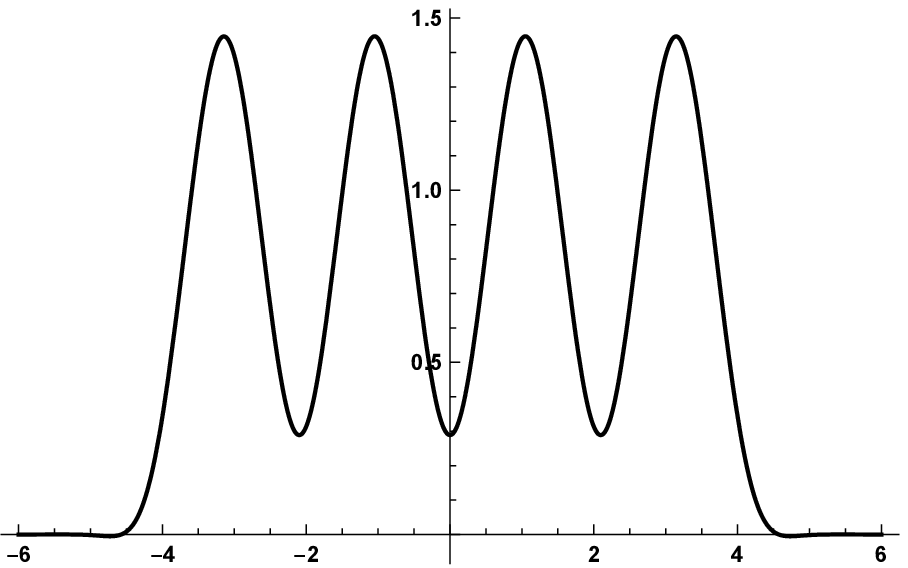} %% mu=-1/4
\includegraphics[height=3.2cm,width=2.6cm]{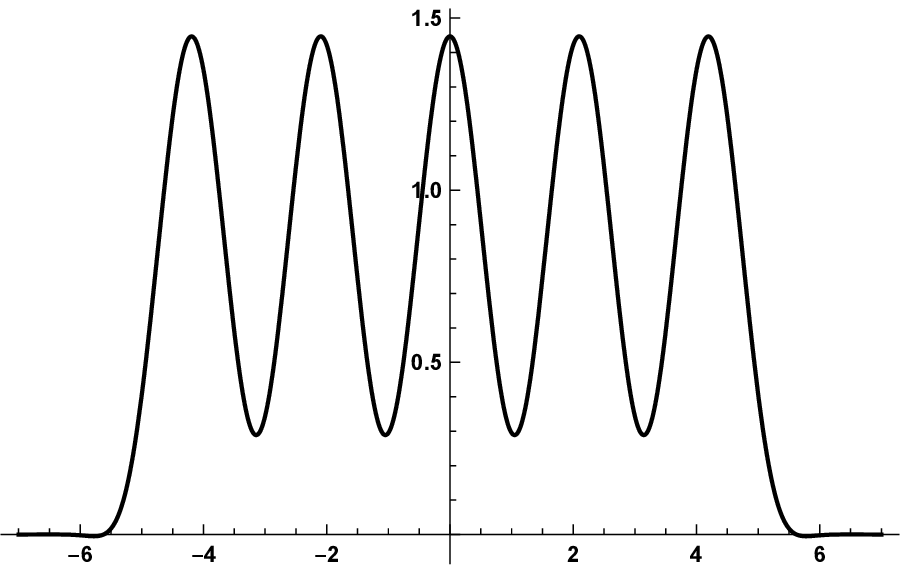} %% mu=-1/4
}
\vskip 0.2cm
\leftline{\hskip 0.3cm\eightpoint{\bf Figure 2.}
$v=\DD u$ for the solutions \EoneNA, \EtwoNA, \EthreeNA, \EfourNB, and \EfiveNC;
for $m=3$ and $\mu=-\sfrac{1}{4}$.}
%%%%%%%%%%%%%%%%%%%%%%%%%%%%%%%%%%%%%%%%%%%%%%%%%%%%%%%%%%%%%

%%%%%%%%%%%%%%%%%%%%%%%%%%%%%%%%%%%%%%%%%%%%%%%%%%%%%%%%%%%%%
\vskip0.6cm
\hbox{\hskip0.4cm
\includegraphics[height=3.2cm,width=3.2cm]{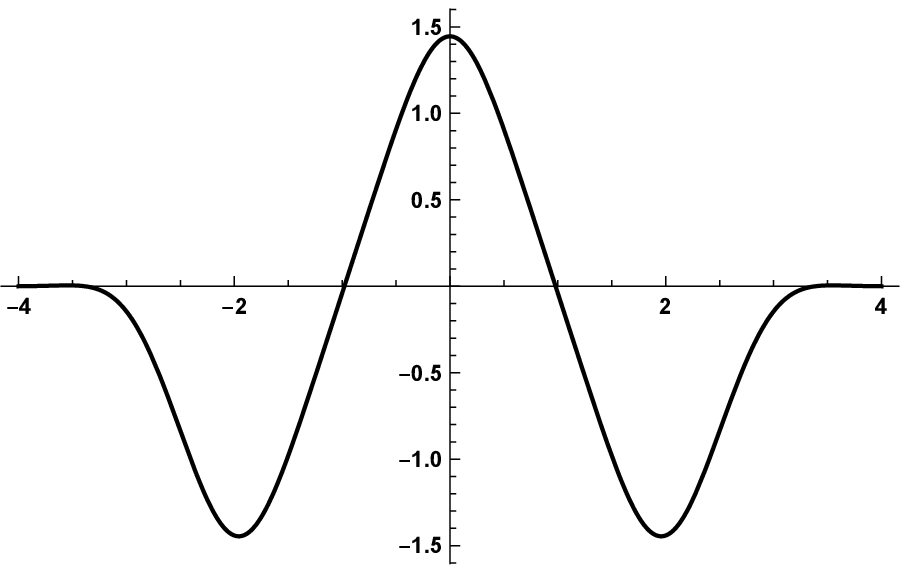} %% mu=-1/4
\includegraphics[height=3.2cm,width=3.2cm]{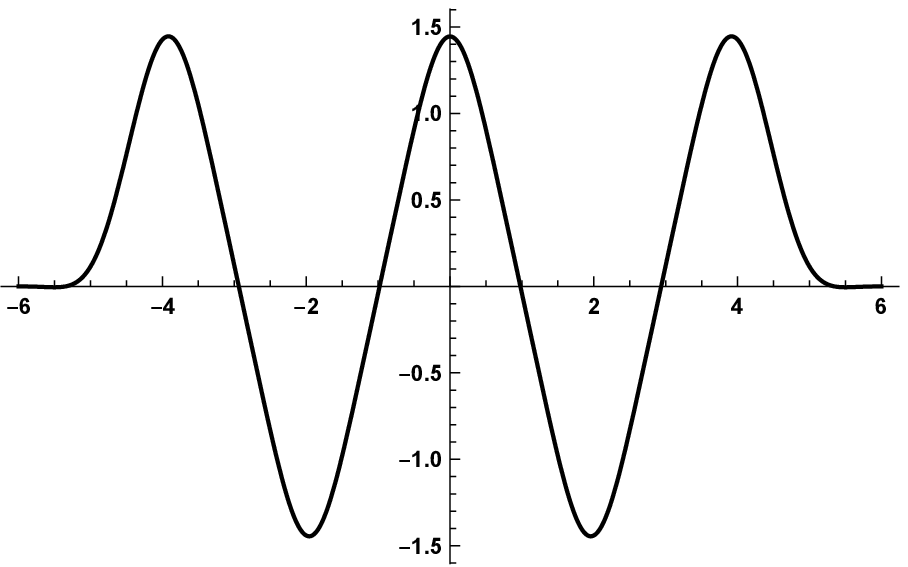} %% mu=-1/4
\includegraphics[height=3.2cm,width=3.2cm]{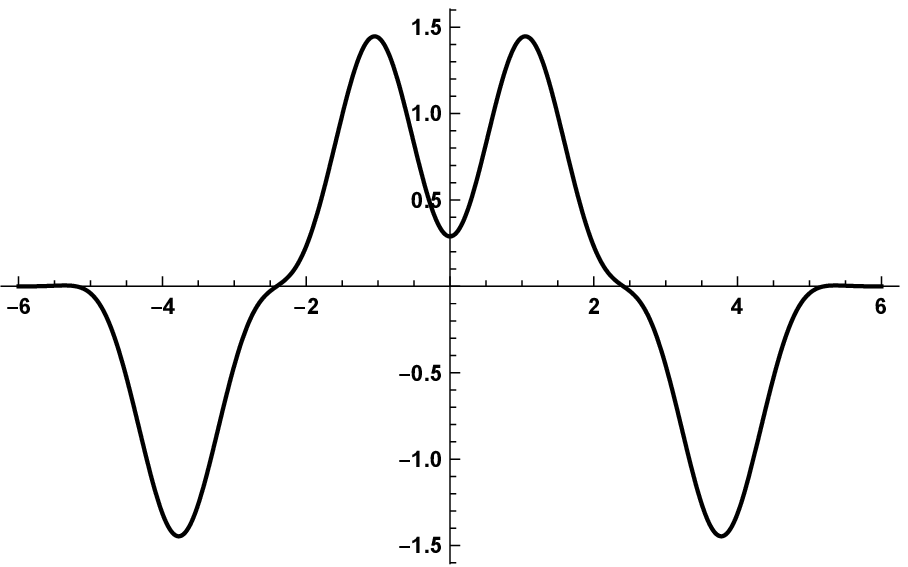} %% mu=-1/4
\includegraphics[height=3.2cm,width=3.2cm]{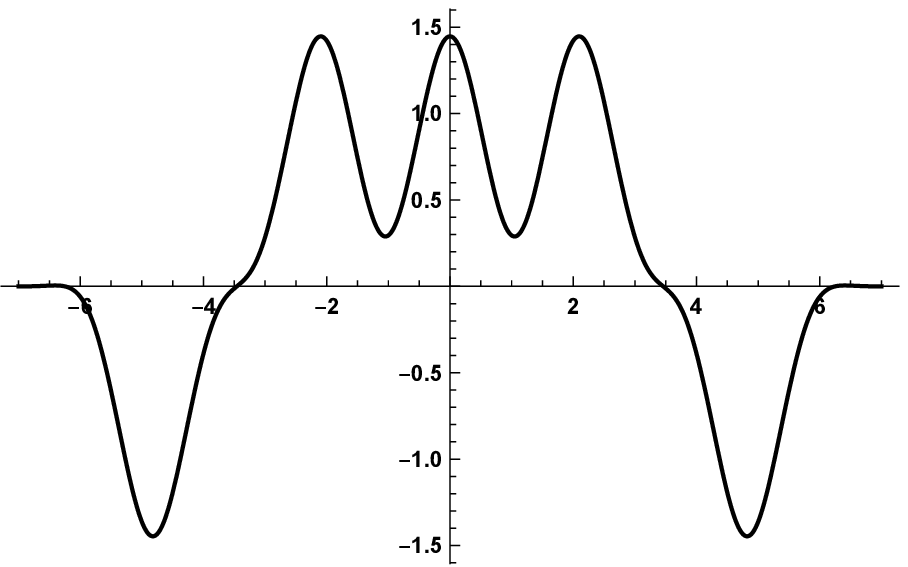} %% mu=-1/4
}
\vskip 0.2cm
\leftline{\hskip 1.0cm\eightpoint{\bf Figure 3.}
$v=\DD u$ for the solutions \EthreeNB, \EfiveNA, \EfourNA, and \EfiveND;
for $m=3$ and $\mu=-\sfrac{1}{4}$.}
%%%%%%%%%%%%%%%%%%%%%%%%%%%%%%%%%%%%%%%%%%%%%%%%%%%%%%%%%%%%%

%%%%%%%%%%%%%%%%%%%%%%%%%%%%%%%%%%%%%%%%%%%%%%%%%%%%%%%%%%%%%
\vskip0.6cm
\hbox{\hskip0.0cm
\includegraphics[height=4.0cm,width=13.0cm]{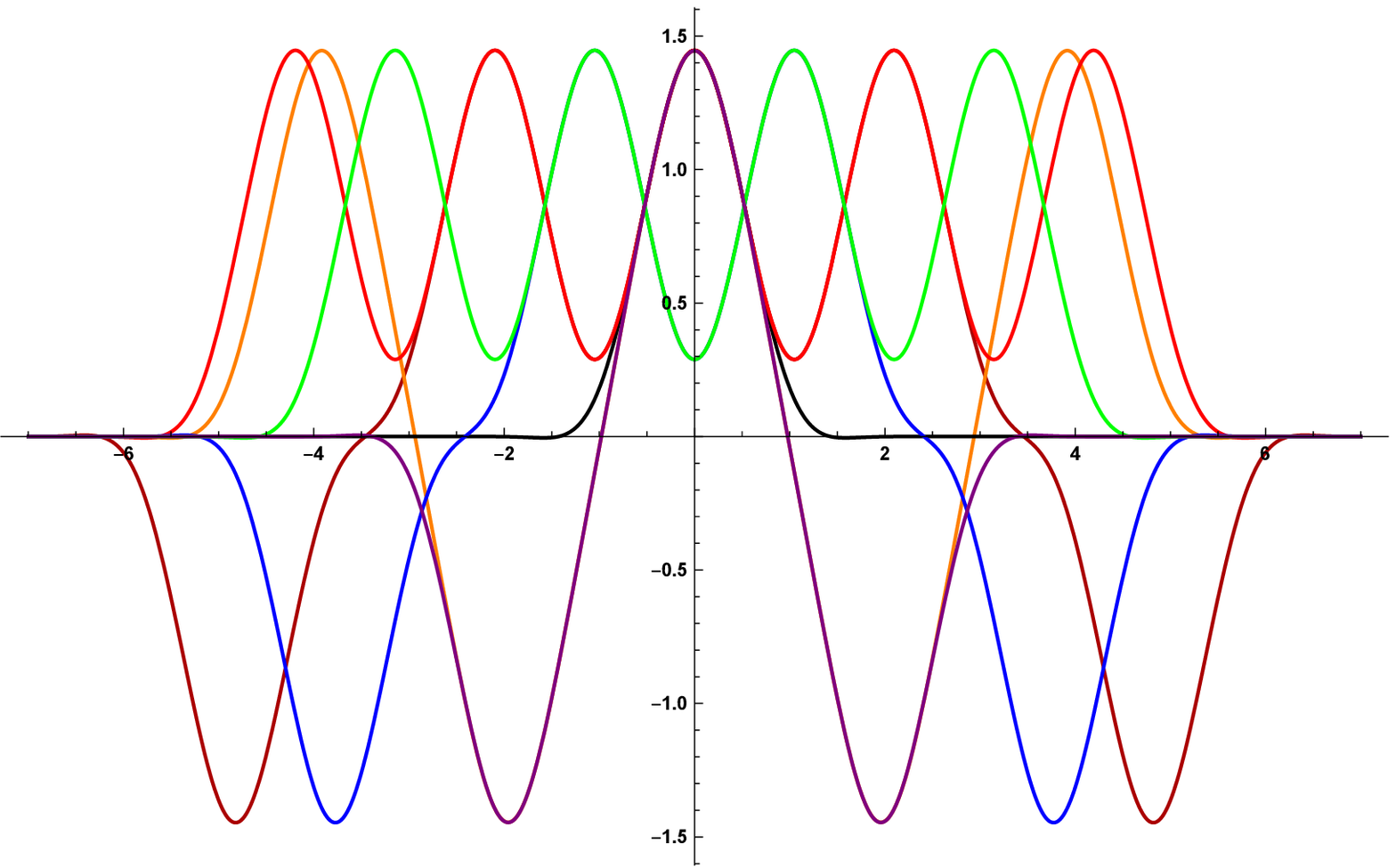}
}\leftline{\hskip 2.0cm\eightpoint{\bf Figure 4.}
Comparison of graphs from Figures 2 and 3.}
%%%%%%%%%%%%%%%%%%%%%%%%%%%%%%%%%%%%%%%%%%%%%%%%%%%%%%%%%%%%%

\subsection Some solutions with $\boldmath\sigma=0$
%%%%%%%%%%%%%%%%%%%%%%%%%%%%%%%%%%%%%%%%%%%%%%%%%%%

%%%%%%%%%%%%%%%%%%%%%%%%%%%%%%%%%%%%%%%%%%%%%%%%%%%%%%%%%%%%%
\vskip0.7cm
\hbox{\hskip0.0cm
\includegraphics[height=3.2cm,width=3.2cm]{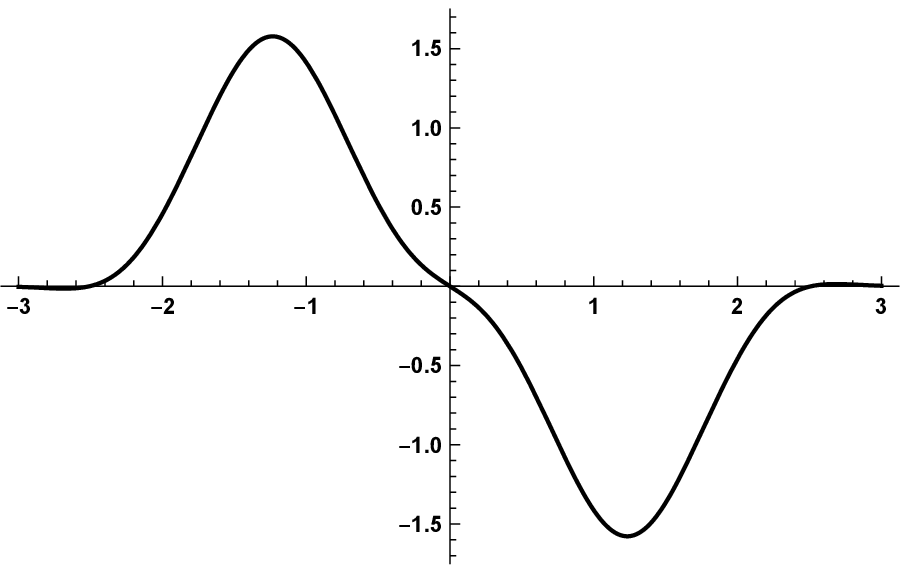} %% mu=-1/2
\includegraphics[height=3.2cm,width=3.2cm]{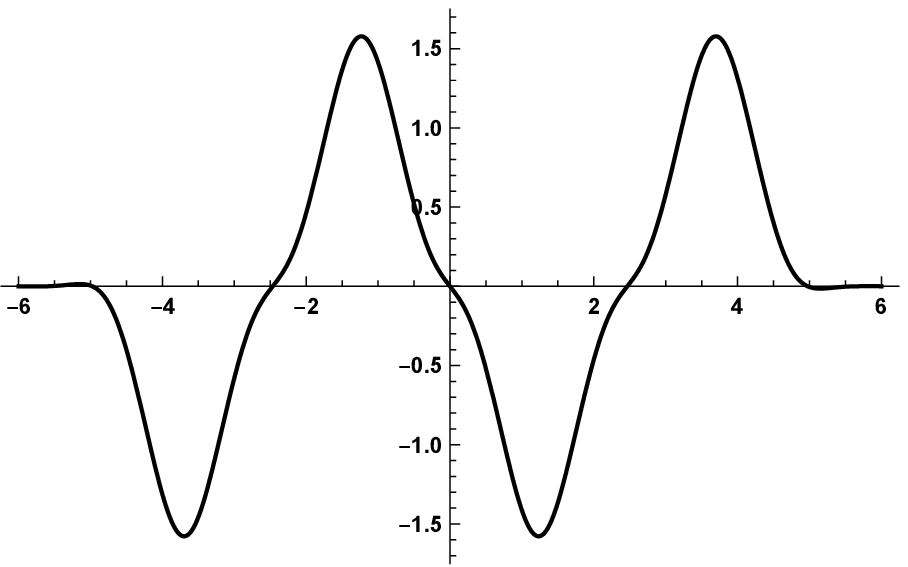} %% mu=-1/2
\includegraphics[height=3.2cm,width=3.2cm]{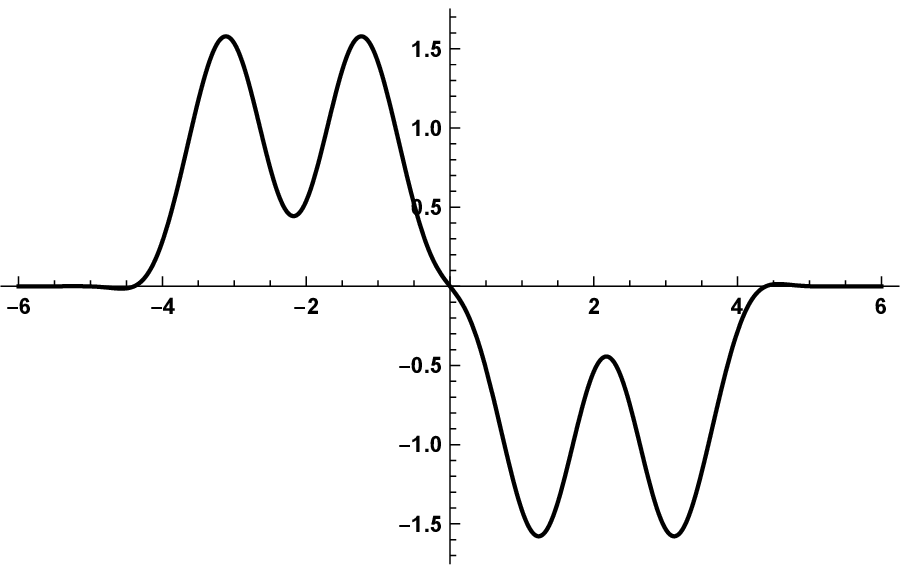} %% mu=-1/2
\includegraphics[height=3.2cm,width=3.2cm]{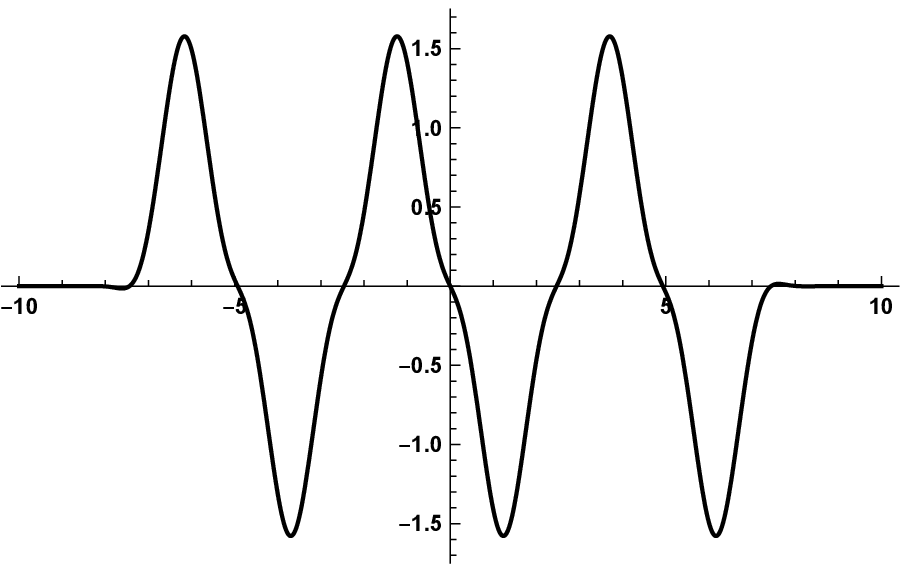} %% mu=-1/2
}
\vskip 0.2cm
\leftline{\hskip 0.1cm\eightpoint{\bf Figure 5.}
$v=\DD u$ for the solutions \OtwoNA, \OfourNA,  \OfourNB, and \OsixNA;
for $m=3$ and $\mu=-\sfrac{1}{2}$.}
%%%%%%%%%%%%%%%%%%%%%%%%%%%%%%%%%%%%%%%%%%%%%%%%%%%%%%%%%%%%%

%%%%%%%%%%%%%%%%%%%%%%%%%%%%%%%%%%%%%%%%%%%%%%%%%%%%%%%%%%%%%
\vskip1.0cm
\hbox{\hskip0.0cm
\includegraphics[height=4.2cm,width=13.0cm]{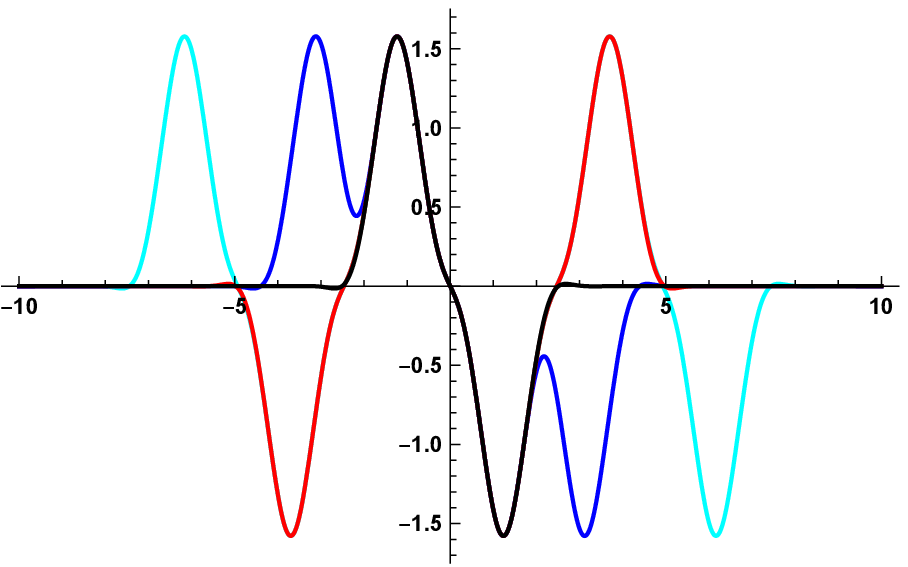}
}\leftline{\hskip 1.5cm\eightpoint{\bf Figure 6.}
Comparison of graphs from Figure 5.}
%%%%%%%%%%%%%%%%%%%%%%%%%%%%%%%%%%%%%%%%%%%%%%%%%%%%%%%%%%%%%

\subsection Limit behavior
%%%%%%%%%%%%%%%%%%%%%%%%%%

%%%%%%%%%%%%%%%%%%%%%%%%%%%%%%%%%%%%%%%%%%%%%%%%%%%%%%%%%%%%%
\vskip0.5cm
\hbox{\hskip1.5cm
\includegraphics[height=3.5cm,width=10cm]{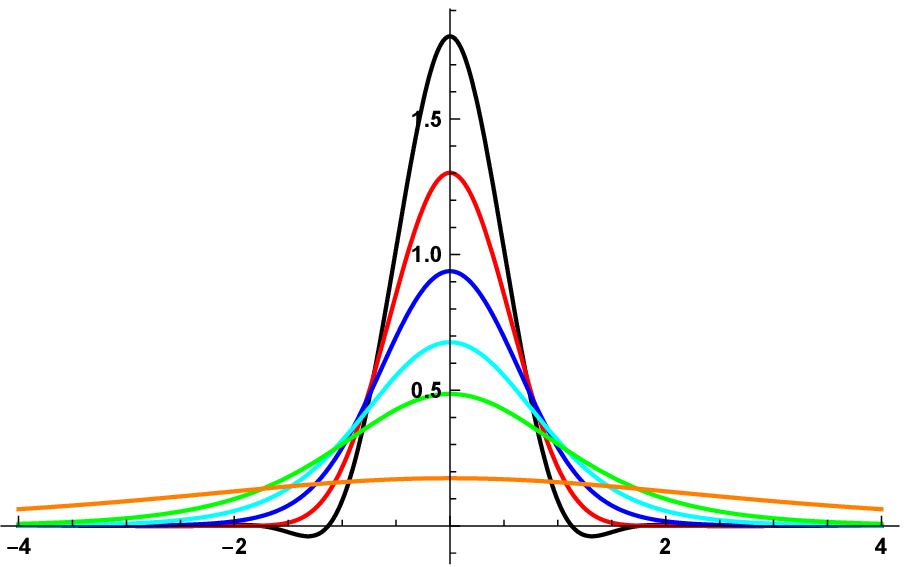}
}
\vskip 0.2cm
\leftline{\hskip 1cm\eightpoint{\bf Figure 7.}
Behavior of $v=\DD u$ as $\mu{\scriptstyle\nearrow\tinyskip} 1$;
$\mu=-\sfrac{63}{64},0,\sfrac{1}{2},\sfrac{3}{4},\sfrac{7}{8},\sfrac{63}{64}$.}
%%%%%%%%%%%%%%%%%%%%%%%%%%%%%%%%%%%%%%%%%%%%%%%%%%%%%%%%%%%%%

%%%%%%%%%%%%%%%%%%%%%%%%%%%%%%%%%%%%%%%%%%%%%%%%%%%%%%%%%%%%%
\vskip1.0cm
\hbox{\hskip0.0cm
\includegraphics[height=4.0cm,width=6cm]{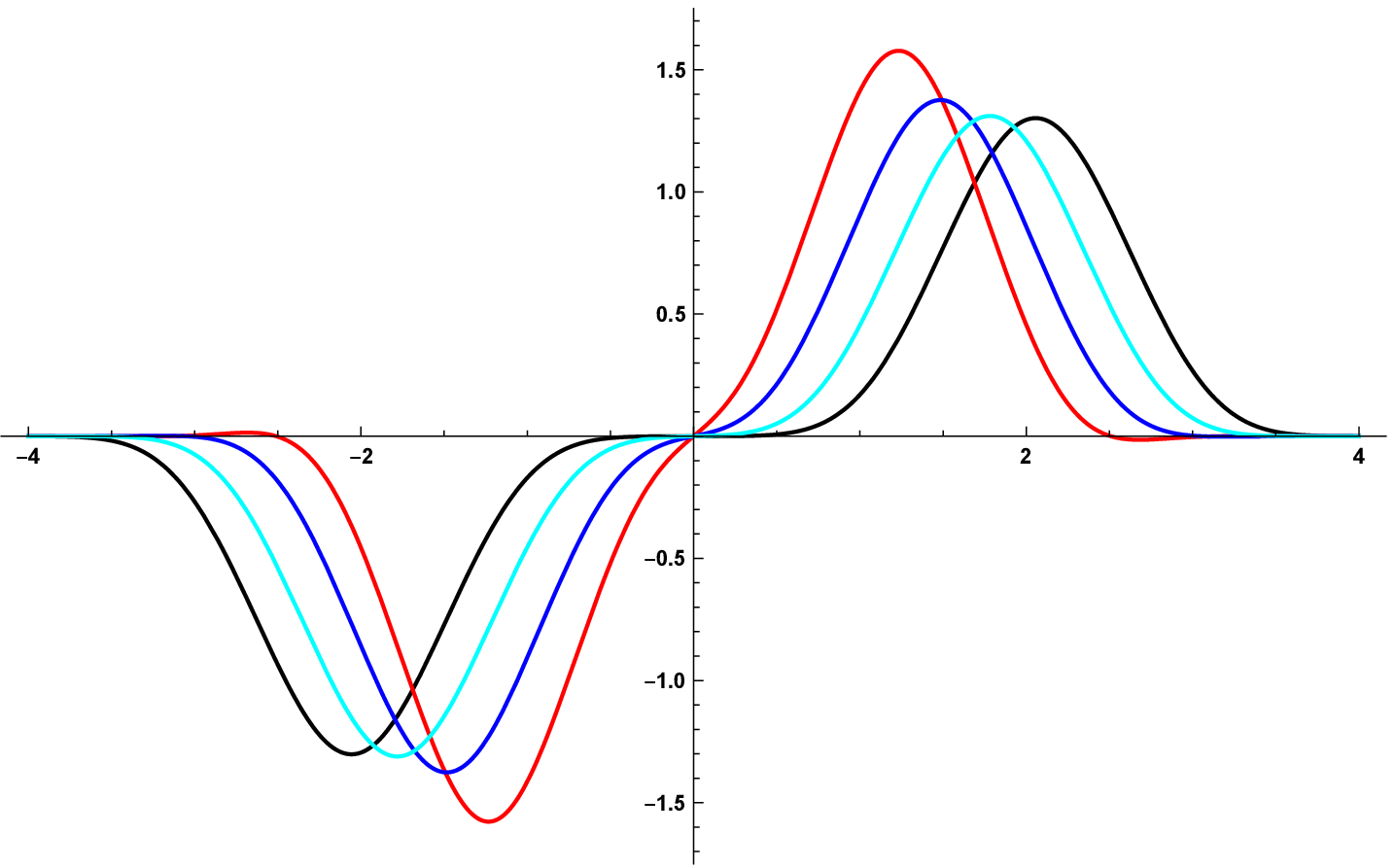}
\includegraphics[height=4.0cm,width=6cm]{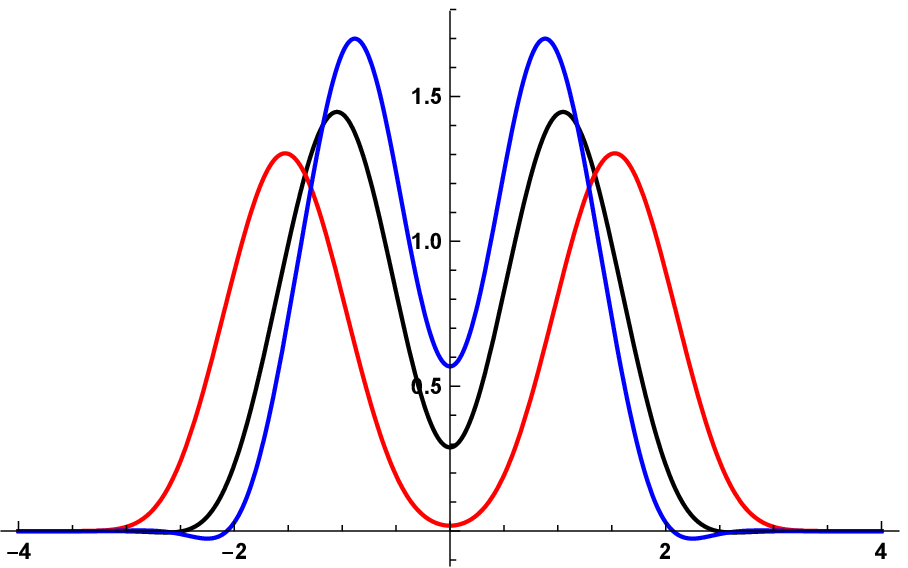}
}
\vskip 0.3cm
\leftline{\hskip 3cm\eightpoint{\bf Figure 8.}
Behavior of $v=\DD u$ as $\mu{\scriptstyle\nearrow\tinyskip} 0$;}
\leftline{\hskip 1cm\eightpoint
$\mu=-\sfrac{1}{2},-\sfrac{1}{8},-\sfrac{1}{64},-\sfrac{1}{1024}$ (left);
$\mu=-\sfrac{1}{2},-\sfrac{1}{4},-\sfrac{1}{256}$ (right).}
%%%%%%%%%%%%%%%%%%%%%%%%%%%%%%%%%%%%%%%%%%%%%%%%%%%%%%%%%%%%%

\bigskip
\references
%%%%%%%%%%%

{\ninepoint

\item{[\rCBH]} H.~Chi, J.~Bell, B.~Hassard,
{\sl Numerical solution of a nonlinear advance-delay differential equation
from nerve conduction theory},
J. Math. Biol. {\bf 24}, 583--601 (1986).
%% The equation models conduction in a myelinated nerve axon ...

\item{[\rAdams]} J.C.~Adams,
{\sl On the expression of the product of any two Legendre's coefficients
by means of a series of Legendre's coefficients},
Proc. R. Soc. London {\bf 27}, 63--71 (1887).

\item{[\rFrWa]} G.~Friesecke, J.A.D.~Wattis,
{\sl Existence theorem for solitary waves on lattices},
Comm. Math. Phys. {\bf 161}, 391--418 (1994).

\item{[\rSmWi]} D.~Smets, M.~Willem,
{\sl Solitary waves with prescribed speed on infinite lattices},
J. Funct. Anal. {\bf 149}, 266--275 (1997).

\item{[\rFPi]} G.~Friesecke, R.L.~Pego,
{\sl Solitary waves on FPU lattices:
I. Qualitative properties, renormalization and continuum limit},
Nonlinearity {\bf 12}, 1601--1627 (1999).

\item{[\rIosK]} G.~Iooss, K.~Kirchg\"assner,
{\sl Travelling waves in a chain of coupled nonlinear oscillators},
Commun. Math. Phys. {\bf 211}, 439--464 (2000).

\item{[\rIos]} G.~Iooss,
{\sl Travelling waves in the Fermi-Pasta Ulam lattice},
Nonlinearity {\bf 13}, 849--866 (2000).

\item{[\rFPii]} G.~Friesecke, R.L.~Pego,
{\sl Solitary waves on FPU lattices: II. Linear implies nonlinear stability},
Nonlinearity {\bf 15} 1343--1359 (2002).

\item{[\rBuDe]} E.~Buksman, J.~De Luca,
{\sl Two-degree-of-freedom Hamiltonian for the time-symmetric two-body problem
of the relativistic action-at-a-distance electrodynamics},
Physical Review E {\bf 67}, 026219 (2003).
%% Wheeler-Feynman electromagnetism theory

\item{[\rFPiii]} G.~Friesecke, R.L.~Pego,
{\sl Solitary waves on FPU lattices: III. Howland-type Floquet theory},
Nonlinearity {\bf 17}, 207--227 (2004).

\item{[\rZg]} P.~Zgliczy\'nski,
{\sl Rigorous numerics for dissipative partial differential equations.
II. Periodic orbit for the Kuramoto-Sivashinsky PDE --- a computer-assisted proof},
Found. Comput. Math. {\bf 4}, 157--185 (2004).

\item{[\rAPankov]} A.~Pankov,
{\sl Traveling Waves And Periodic Oscillations in Fermi-Pasta-Ulam Lattices},
Imperial College Press (2005).

\item{[\rJaSii]} G.~James, Y.~Sire,
{\sl Travelling breathers with exponentially small tails in a chain of nonlinear oscillators},
Comm. Math. Phys. {\bf 257}, 51--85 (2005).

\item{[\rGGKKMO]} M.~Gameiro, T.~Gedeon, W.~Kalies, H.~Kokubu, K.~Mischaikow, H.~Oka,
{\sl Topological horseshoes of traveling waves for a fast-slow predator-prey system},
J. Dynam. Differ. Equations {\bf 19} 23--654 (2007).

\item{[\rKaAl]} A.~Kaddar, H.T.~Alaoui,
{\sl Fluctuations in a mixed IS-LM business cycle model},
Electron. J. Differential Equations {\bf 134}, 1--9 (2008).

\item{[\rMiNa]} T.~Minamoto, M.T.~Nakao,
{\sl A numerical verification method for a periodic solution of a delay differential equation},
J. Comput. Appl. Math. {\bf 235}, 870--878, (2010).

\item{[\rHeRa]} M.~Herrmann, J.D.M.~Rademacher,
{\sl Heteroclinic travelling waves in convex FPU-type chains},
SIAM J. Math. Anal. {\bf 42}, 1483--1504 (2010).

\item{[\rAKi]} G.~Arioli, H.~Koch,
{\sl Integration of dissipative PDEs: a case study},
SIAM J. Appl. Dyn. Syst. {\bf 9} 1119--1133 (2010).

\item{[\rPaRo]} A.~Pankov, V.M.~Rothos,
{\sl Traveling waves in Fermi-Pasta-Ulam lattices with saturable nonlinearities},
Discrete Contin. Dynam. Systems A {\bf 30}, 835--849 (2011).

\item{[\rHoWa]} A.~Hoffman, C.E.~Wayne,
{\sl A simple proof of the stability of solitary waves
in the Fermi--Pasta--Ulam model near the KdV limit},
Fields Inst. Comm. {\bf 64}, 185--192 (2013).

\item{[\rAKii]} G.~Arioli, H.~Koch,
{\sl Existence and stability of traveling pulse solutions
of the FitzHugh--Nagumo equation},
Nonlinear Anal. {\bf 113}, 51--70 (2015).

\item{[\rCZg]} A.~Czechowski, P.~Zgliczy\'nski,
{\sl Existence of periodic solutions of the FitzHugh--Nagumo Equations
for an Explicit Range of the Small Parameter},
SIAM J. Appl. Dyn. Syst. {\bf 15}, 1615--1655 (2016).

\item{[\rSZg]} R.~Szczelina, P.~Zgliczy\'nski,
{\sl Algorithm for rigorous integration of delay differential equations
and the computer-assisted proof of periodic orbits in the Mackey-Glass equation},
Found. Comput. Math. (2017).

\item{[\rKLM]} J.~Jaquette, J.-P.~Lessard, K.~Mischaikow,
{\sl Stability and uniqueness of slowly oscillating periodic solutions to Wright's equation},
J. Differ. Equations {\bf 263}, 7263--7286 (2017).

\item{[\rKiLe]} G.~Kiss, J.-P.~Lessard,
{\sl Rapidly and slowly oscillating periodic oscillations of a delayed van der Pol oscillator},
J. Dynam. Differ. Equations {\bf 29}, 1233--1257 (2017).

\item{[\rMatsue]} K.~Matsue,
{\sl Rigorous numerics for fast-slow systems with one-dimensional slow variable:
topological shadowing approach},
Topol. Methods Nonlinear Anal. {\bf 50}, 357--468 (2017).

\item{[\rAKiv]} G.~Arioli, H.~Koch,
{\sl Non-radial solutions for some semilinear elliptic equations on the disk},
Nonlin. Anal. TMA {\bf 179}, 294--308 (2019).

\item{[\rAKv]} G.~Arioli, H.~Koch,
{\sl Some breathers and multi-breathers for FPU-type chains},
\pdfclink{0 0 1}{Preprint 2018}{http://web.ma.utexas.edu/users/koch/papers/breathers/},
to appear in Commun. Math. Phys.

\item{[\rWZg]} D.~Wilczak, P.~Zgliczy\'nski,
{\sl Symbolic dynamics for Kuramoto-Sivashinsky PDE on the line --- a computer-assisted proof},
Preprint 2017, arXiv:1710.00329.

%%%%%%%%% COMP

\item{[\rFiles]} G.~Arioli, H.~Koch.
The source code for our programs, and data files, are available at\hfill\break
\pdfclink{0 0 1}{{\tt web.ma.utexas.edu/users/koch/papers/twaves/}}
{http://web.ma.utexas.edu/users/koch/papers/twaves/}

\item{[\rAda]} Ada Reference Manual, ISO/IEC 8652:2012(E),
available e.g. at\hfil\break
\pdfclink{0 0 1}{{\tt www.ada-auth.org/arm.html}}
{http://www.ada-auth.org/arm.html}

\item{[\rGnat]}
A free-software compiler for the Ada programming language,
which is part of the GNU Compiler Collection; see
\pdfclink{0 0 1}{{\tt gnu.org/software/gnat/}}{http://gnu.org/software/gnat/}

\item{[\rIEEE]} The Institute of Electrical and Electronics Engineers, Inc.,
{\sl IEEE Standard for Binary Float\-ing--Point Arithmetic},
ANSI/IEEE Std 754--2008.

\item{[\rMPFR]} The MPFR library for multiple-precision floating-point computations
with correct rounding; see
\pdfclink{0 0 1}{{\tt www.mpfr.org/}}{http://www.mpfr.org/}

}

\bye